\documentclass[12pt,reqno]{amsart}
\usepackage{amsmath}
\usepackage{amssymb}
\usepackage{amsthm}
\usepackage{latexsym}
\usepackage{amssymb}

\newtheorem{theorem}{Theorem}[section]
\newtheorem{lemma}{Lemma}[section]
\newtheorem{corollary}{Corollary}[section]
\newtheorem{remark}{Remark}[section]
\newtheorem{definition}{Definition}[section]
\newtheorem{proposition}{Proposition}[section]
\newtheorem{example}{Example}[section]
\newtheorem{assumption}{Assumption}[section]
\numberwithin{equation}{section}
\newcommand{\bth}{\begin{theorem}}
\newcommand{\ethe}{\end{theorem}}

\newcommand{\bre}{\begin{remark}}
\newcommand{\ere}{\end{remark}}

\newcommand{\ble}{\begin{lemma}}
\newcommand{\ele}{\end{lemma}}
\newcommand{\bde}{\begin{definition}}
\newcommand{\ede}{\end{definition}}

\newcommand{\bco}{\begin{corollary}}
\newcommand{\eco}{\end{corollary}}

\newcommand{\bpr}{\begin{proposition}}
\newcommand{\epr}{\end{proposition}}

\newcommand{\bexer}{\begin{exercise}}
\newcommand{\eexer}{\end{exercise}}
\newcommand{\breh}{\begin{hint}}
\newcommand{\ereh}{\end{hint}}

\newcommand{\halmos}{\hfill \qed}

\newcommand{\bexam}{\begin{example}}
\newcommand{\eexam}{\end{example}}

\newcommand{\pr} {{\bf Proof.}}

\newcommand{\bfi}{\begin{fig}}
\newcommand{\efi}{\end{fig}}

\newcommand{\beao}{\begin{eqnarray*}}
\newcommand{\eeao}{\end{eqnarray*}\noindent}

\newcommand{\beam}{\begin{eqnarray}}
\newcommand{\eeam}{\end{eqnarray}\noindent}

\newcommand{\PP}{\mathbf{P}}
\newcommand{\nto}{n\to\infty}

\newcommand{\xto}{x\to\infty}

\newcommand{\bH}{\overline{H}}
\newcommand{\bF}{\overline{F}}
\newcommand{\bB}{\overline{B}}
\newcommand{\bG}{\overline{G}}

\newcommand{\bbr}{{\mathbb R}}

\newcommand{\bbb}{{\mathbb B}}

\newcommand{\bbn}{{\mathbb N}}

\newcommand{\yto}{y\to\infty}

\newcommand{\vep}{\varepsilon}

\begin{document}
\title{A new approach in two-dimensional heavy-tailed distributions}

\author[ D.G. Konstantinides, C. D. Passalidis ]{ Dimitrios G. Konstantinides, Charalampos  D. Passalidis}

\address{Dept. of Statistics and Actuarial-Financial Mathematics,
University of the Aegean,
Karlovassi, GR-83 200 Samos, Greece}
\email{konstant@aegean.gr,\;sasm23002@sas.aegean.gr.}

\date{{\small \today}}

\begin{abstract}
We consider a new approach in the definition of two-dimensional heavy-tailed distributions. Namely, we introduce the classes of two-dimensional long-tailed, of two-dimensional dominatedly varying and of two-dimensional consistently varying distributions. Next, we define the closure property with respect to two-dimensional convolution and to joint max-sum equivalence in order to study if they are satisfied by these classes. Further we examine the joint behavior of two random sums, under generalized tail asymptotic independence. Afterward we study the closure property under scalar product and two dimensional product convolution and by these results we extended our main result in the case of jointly randomly weighted sums. Our results contained some applications where we establish the asymptotic expression of the ruin probability in a two-dimensional discrete-time risk model. 
\end{abstract}

\maketitle
\vspace{3mm}


\section{Introduction} \label{sec.KP.1}

\subsection{Preliminaries} \label{subsec.KP.1.1}

The heavy-tailed distributions describe precisely complicated situations. One of most important application is related to the risk theory in actuarial science. Although several one-dimensional problems remain still open, the multidimensional case meets popularity from both theoretical and practical aspect. Especial, with respect to practical point of view, the modern insurance industry does not operate with a single portfolio. 

On this line there are some recent papers, as for example \cite{hu:jiang:2013},  \cite{konstantinides:li:2016}, \cite{yang:su:2023}. On this direction, we introduce some two-dimensional distribution classes, with heavy tails, that are convenient to calculations and permit direct and consistent generalization of the one-dimensional concepts.

In subsection 1.2, we remind some basic definitions, for one-dimensional heavy-tailed distributions, for easy comparison with the two-dimensional ones. In section 2, we introduce the closure property with respect to the two-dimensional convolution and the two-dimensional max-sum equivalence. Next, we present some results on these classes of distributions. In section 3, we estimate the joint asymptotic behavior of two random sums, under a dependence structure that generalizes the tail asymptotic independence, and we establish an asymptotic expression for the ruin probabilities, in a discrete-time two-dimensional risk model without stochastic discount factors. Furthermore in section 5 we study the closure property of some of new classes with respect to scalar product, and in section 6 we extended some of our results in section 4, in the case wich we have a common discount factor for the two portfolios. Last but not least we limited ourselves in the non-negative case and we study the closure property of new classes with respect to product convolution in two dimensions, and some previous results are extended.

Before passing to the next subsection we give some notations, that we need for the rest of the paper. We denote by $\bF:=1-F$ the distribution tail, hence $\bF(x)=\PP[X>x]$ and holds $\bF(x)>0$
for any $x \geq 0$, except it is referred differently. For two positive functions $f(x)$ and $g(x)$, the asymptotic relation $f(x)=o[g(x)]$, as $\xto$ means
\beao
\lim_{\xto}\dfrac{f(x)}{g(x)} = 0\,,
\eeao
the asymptotic relation $f(x)=O[g(x)]$, as $\xto$ holds if
\beao
\limsup_{\xto} \dfrac{f(x)}{g(x)} < \infty\,.
\eeao
and the asymptotic relation $f(x)\asymp g(x)$, as $\xto$ if  both $f(x)=O[g(x)]$ and 
$g(x)=O[f(x)]$. Similarly, for the bi-variate functions $f(x,\,y)$, $g(x,\,y)$ hold the corresponding asymptotic relations with $\min\{x,\,y\} \to \infty$, as for example $f(x,\,y) = o[g(x,\,y)]$ if it holds
\beao
\lim_{x \wedge y \to \infty} \dfrac{f(x,\,y) }{g(x,\,y) }=0\,.
\eeao
For a real number $x,y$, we denote $x^{+}:=\max\{x,0\}$, $x\wedge y:=\min\{x,y\}$,  $x\vee y:=\max\{x,y\}$. With bold letters, we denote vectors and further for the unit and zero vectors we write ${\bf 1}$ and ${\bf 0}$, respectively.

\subsection{Uni-variate heavy-tailed distributions}

The following properties are to be extended in two dimensions.
\begin{enumerate}
\item
For two random variables $X_1$, $X_2$ with distributions $F_1$, $F_2$ respectively, the distribution of the sum is defined by $F_{X_1+X_2}(x)=\PP[X_1 + X_2 \leq x]$ with tail $\overline{F_{X_1+X_2}}(x)=\PP[X_1+X_2 >x]$. If $X_1$, $X_2$ are independent, we write $F_1*F_2$ instead of $F_{X_1+X_2}$. 
\item
We say that the random variables  $X_1$, $X_2$ or their distributions $F_1$, $F_2$ are max-sum equivalent if $\overline{F_1*F_2}(x) \sim \bF_1(x)+\bF_2(x)$, as $\xto$.
\end{enumerate}

Now we consider some classes of heavy-tailed distributions. We say that a distribution $F$ is heavy-tailed, and we write $F \in \mathcal{K}$, if holds
\beao
\int_{-\infty}^{\infty} e^{\vep\,x}\,F(dx) = \infty\,,
\eeao
for any $\vep>0$. A large enough class of heavy-tailed distributions is the class of long tails, denoted by $\mathcal{L}$. We have $F \in \mathcal{L}$ if holds
\beao
\lim_{\xto} \dfrac{\bF(x-a)}{\bF(x)}=1\,,
\eeao
for any (or equivalently, for some) $a> 0$. It is well-known that if $F \in \mathcal{L}$, then there exists a function $a\;:\;[0,\,\infty) \longrightarrow [0,\,\infty)$, such that $a(x) \rightarrow \infty$, $\bF(x\pm a(x)) \sim \bF(x)$, as $\xto$. This kind of function $a(x)$ is called insensitivity function for $F$, see further in \cite{cline:samorodnitsky:1994}, \cite{foss:korshunov:zachary:2013} or \cite{Konstantinides:2018}.

A little smaller class than $\mathcal{L}$ is the class of subexponential distributions, introduced in \cite{chistyakov:1964}. We say that a distribution $F$ with support the interval $[0,\,\infty)$ belongs to the class of subexponential distributions, symbolically $F \in \mathcal{S}$ if it holds
\beao
\lim_{\xto} \dfrac{\overline{F^{n*}}(x)}{\bF(x)}=n\,,
\eeao
for any $n \in \bbn$, where $F^{n*}$ represents the $n$-th order convolution power for $F$. The class $\mathcal{S}$ has found several applications in the risk models, as for example in \cite{li:tang:wu:2010}, \cite{geng:liu:wang:2023}, \cite{ji:wang:yan:cheng:2023}.

We say that the distribution $F$ belongs to the class of the dominatedly varying distributions, symbolically $F \in \mathcal{D}$, if holds
\beao
\limsup_{\xto} \dfrac{\bF(b\,x)}{\bF(x)} < \infty\,,
\eeao
for some (or equivalently, for all) $b\in (0,\,1)$. Is well known that $\mathcal{D}\cap \mathcal{L}=\mathcal{D}\cap \mathcal{S} \subset \mathcal{K}$, see \cite[Th.1]{goldie:1978}.

Further, a smaller class of heavy-tailed distributions represents the class of consistently varying distributions, symbolically $F\in \mathcal{C}$. We say that $F\in \mathcal{C}$, if holds
\beao
\lim_{y\uparrow 1}\limsup_{\xto} \dfrac{\bF(y\,x)}{\bF(x)} =1\,,
\eeao 
or equivalently
\beao
\lim_{y\downarrow 1}\liminf_{\xto} \dfrac{\bF(y\,x)}{\bF(x)} =1\,.
\eeao

Finally, we say that a distribution $F$ belongs to the class of regularly varying distributions, with index $\alpha>0$,  symbolically $F \in \mathcal{R}_{-\alpha}$ if holds
\beao
\lim_{\xto} \dfrac{\bF(t\,x)}{\bF(x)} =t^{-\alpha}\,,
\eeao 
for any $t>0$.

For these classes we obtain the following inclusions (see \cite{bingham:goldie:teugels:1987})
\beao
\mathcal{R}:=\bigcup_{\alpha \geq 0} \mathcal{R}_{-\alpha} \subsetneq \mathcal{C} \subsetneq \mathcal{D}\cap \mathcal{L} \subsetneq \mathcal{S} \subsetneq \mathcal{L} \subsetneq \mathcal{K}\,,
\eeao 
where $\mathcal{R}_0$ is the class of slowly varying distributions. We can find numerous classes of heavy-tailed distributions, however we mentioned the most popular in the literature. In this paper we extend into two dimensions the classes $ \mathcal{C}$, $ \mathcal{D}$ and $ \mathcal{L}$.

In \cite{cai:tang:2004} we find the following results.

\bpr \label{pr.KP.1}
If $F_1 \in \mathcal{D}$ and $F_2 \in \mathcal{D}$ are distributions with support the interval $[0,\,\infty)$, then $F_{X_1+X_2} \in \mathcal{D}$.
\epr

In Proposition \ref{pr.KP.1} we find that for non-negative random variables, the class $\mathcal{D}$ satisfies the closure property with respect to sum. As was mentioned in \cite{cai:tang:2004}, the class $\mathcal{D}$ does NOT satisfy the max-sum equivalence, as it follows from the fact that $\mathcal{D} \not\subset \mathcal{S}$ and $\mathcal{S} \not\subset \mathcal{D}$, therefore the relation $\overline{F^{2*}}(x) \sim 2\,\bF(x)$, as $\xto$, does NOT hold for $F\in \mathcal{D} \setminus \mathcal{S}$. In opposite to the dominated variation, the class of the consistently varying  distributions satisfy both these properties.

\bpr \label{pr.KP.2}
If $F_1 \in \mathcal{C}$ and $F_2 \in \mathcal{C}$ are distributions with support the interval $[0,\,\infty)$, then it holds $F_1 * F_2 \in \mathcal{C}$ and $\overline{F_1* F_2}(x) \sim \bF_1(x) + \bF_2(x)$, as $\xto$.
\epr

\section{Two-dimensional heavy tails} \label{sec.KP.2}

The reason why the multivariate distributions have been so popular is their ability to describe better multidimensional phenomena. This happens because of the interdependence among the components of the random vectors, the affect significantly on the final outcome.

The first heavy-tailed distributions class that was extended to multidimensional frame is the regular variation. We say that the random vector ${\bf X}=(X_1,\,\ldots,\,X_d)$ represents a multivariate regularly varying vector with index $\alpha$ and Radon measure $\nu$, symbolically ${\bf X} \in MRV(\alpha,\,F,\,\nu)$ if holds
\beao
\lim_{\xto} \dfrac 1{\bF(x)}\PP\left[ \dfrac {{\bf X}}x \in \bbb \right]=\nu(\bbb)\,,
\eeao
for any $\nu$-continuous Borel set $\bbb \subset [0,\,\infty]^d\setminus \{ {\bf 0}\} $, with $F \in \mathcal{R}_{-\alpha}$ and the measure $\nu$ is homogeneous, namely holds $\nu (\lambda\,\bbb) = \lambda^{-\alpha}\,\nu(\bbb)$, for any  $\lambda > 0$.

The frame of multivariate regular variation was introduced in \cite{dehaan:resnick:1982}. Under this definition, the multivariate regular variation was used in the study of several issues in multivariate risk models and in risk management, as for example in \cite{li:2016}, \cite{tang:yang:2019}, \cite{yang:su:2023}.

Although this kind of extension to multidimensional setup is well-established, it does not happen to other multidimensional distribution classes. Most of the extensions cover the multivariate subexponential distribution class and the multivariate long tailed distribution class.

Initially in \cite{cline:resnick:1992} was introduced these two distribution classes, as essential extension of the multivariate regular variation, namely using vague convergence and point processes. Later, in \cite{omey:2006} appear three different formulations for the multivariate subexponentiality and the multivariate long-tailedness. The formulations, that are close to our definitions, are given in classes $\mathcal{S}(\bbr^d)$ and $\mathcal{L}(\bbr^d)$.   
We say that the multivariate distribution $F$ belongs to class $\mathcal{S}(\bbr^d)$, if it holds
\beao
\lim_{\xto} \dfrac{\overline{F^{2*}}({\bf t}\,x)}{\bF({\bf t}\,x)}=2\,,
\eeao
for any ${\bf x} > {\bf 0}$, with $\min_{1\leq i \leq d} \{t_i\} < \infty$, and that the multivariate distribution $F$ belongs to class $\mathcal{L}(\bbr^d)$, if it holds
\beao
\lim_{\xto} \dfrac{\bF({\bf t}\,x-{\bf a})}{\bF({\bf t}\,x)}=1\,,
\eeao
for any ${\bf a} \geq  {\bf 0}$ and for any ${\bf t} > {\bf 0}$, with $\min_{1\leq i \leq d} \{t_i\} < \infty$.

This approach was used to study the asymptotic behavior of the tail of randomly stopped sum of random vectors, namely $S_{N}=\sum_{i=1}^N {\bf X}_i$, where $N$ is a discrete random variable with support $\bbn_0=\bbn \cup \{0\}$ and the $ {\bf X}_i$ are independent, identically distributed random vectors with multivariate distribution $F$. For applications of this class, see in \cite{omey:mallor:santos:2006}.

Finally, another formulation of multivariate subexponential distributions was provided in \cite{samorodnitsky:sun:2016}, which represents the only approach with results for the ruin probability in a multivariate continuous-time risk model.

In the present paper we confine ourselves in the two dimensions and we stay close to the formulation in \cite{omey:2006}, however we keep two important differences.

At first, we follow a direct approach to the uni-variate distribution classes definitions.

At second, in the case of $d=2$ the formulation in \cite{omey:2006}, and in the definition of multivariate regular variation,  the convention ${\bf F}(x,\,y)=\PP[X\leq x,\,Y\leq y]$ is adopted and the distribution tail $1-{\bf F}(x,\,y)$, denoted $\overline{{\bf F}}(x,\,y)$, is applied on the event $\{X>x\} \cup \{Y>y\}$. We consider only the case, in which there exist excesses of both random variables $\{X>x\} \cap \{Y>y\}$, namely we define by ${\bf \bF_1}(x,\,y):=\PP[X>x,\,Y>y]$, as the distribution tail of ${\bf F}$, with notation ${\bf \bF_b}(x,\,y):=\PP[X>b_1\,x,\,Y>b_2\,y]$, for all ${\bf b}=(b_1,\,b_2) \in (0,\,\infty)^2$. The choice of such a definition is due to both the consistency with the uni-variate case and the easiness in asymptotic calculation of the joint tail of random sums as well. We intent that our approach becomes more consistent with the ruin of all portfolios, which represents the worst event which can happen for an insurance company with multiple businesses. 

Next, we introduce the bi-variate heavy-tailed distribution class. From now and further by the notation ${\bf a}=(a_1,\,a_2)> (0,\,0)$ we mean that $(a_1,\,a_2) \in [0,\,\infty)^2 \setminus \{{\bf 0}\}$, except it is referred differently.

\bde \label{def.KP.1}
We say that the random pair (X,\,Y) with marginal distributions $F$, $G$ belongs to the bi-variate long-tailed distributions, symbolically $(F,\,G) \in \mathcal{L}^{(2)}$, if hold
\begin{enumerate}
\item
$F \in \mathcal{L}$ and $G \in \mathcal{L}$.
\item
It holds
\beao
\lim_{x\wedge y \to \infty} \dfrac{{\bf \bF_1}(x-a_1,\,y-a_2)}{{\bf \bF_1}(x,\,y)}=\lim_{x\wedge y \to \infty} \dfrac{\PP[X>x-a_1,\,Y>y-a_2]}{\PP[X>x,\,Y>y]} =1\,,
\eeao
for some, or equivalently for any, ${\bf a}=(a_1,\,a_2) > (0,0)$, with $a_1$ not necessarily equal to $a_2$.
\end{enumerate}
\ede 

\bre \label{rem.KP.2.0} 
From the previous definition we wonder if by the two-dimensional property of class $\mathcal{L}^{(2)}$ follows directly the inclusion $F,\,G \in \mathcal{L}$, while the opposite question is to be replied by Example \ref{exam.KP.2.2}. The answer to this question is no, because it holds for any, or equivalently for some, $(a_1,\,a_2) >(0,\,0)$, as follows from Definition \ref{def.KP.1}.

Let $F \in  \mathcal{L}$ be a distribution and $G$ be another distribution, not necessarily from class $\mathcal{L}$. We assume that the two distributions stem from the independent random variables $X$ and $Y$, thus if $a_1>0$ and $a_2=0$ we find that
\beao
\lim_{x\wedge y \to \infty} \dfrac{\PP[X>x-a_1\,, \;Y>y]}{\PP[X>x\,, \;Y>y]}=\lim_{x\wedge y \to \infty} \dfrac{\bF(x-a_1)\,\bG(y)}{\bF(x)\,\bG(y)}=1\,,
\eeao
however, if it holds $G \notin \mathcal{L}$, then we have not this pair in the class $\mathcal{L}^{(2)}$.

The reason, why we require that the marginals belong to class  $\mathcal{L}$, is to secure some  two-dimensional closure properties, that could fail if the  $\mathcal{L}$ condition is missing. 
\ere

From Definition \ref{def.KP.1}  we obtain that, if $(F,\,G) \in \mathcal{L}^{(2)}$, then for any $(A_1,\,A_2)>(0,\,0)$ holds
\beam \label{eq.KP.2.8}
\sup_{|a_1|<A_1,\,|a_2|<A_2}\left| \PP[X>x-a_1,\,Y >y-a_2]-\PP[X>x,\,Y>y]\right|=o\left( \PP[X>x,\,Y>y]\right),
\eeam 
as $x\wedge y \to \infty $, which follows from the uniformity of the convergence
\beao
\lim_{x\wedge y \to \infty }\dfrac{ \PP[X>x-a_1\,,\;Y >y-a_2]}{\PP[X>x\,,\;Y>y]}=1\,,
\eeao
over the parallelogram $[-A_1\,,\;A_1]\times [-A_2\,,\;A_2] $. Indeed, for $-A_1 \leq a_1 \leq A_1$ and $-A_2 \leq a_2 \leq A_2$ we obtain
$x-A_1 \leq x+a_1 \leq x+A_1$ and $y-A_2 \leq y+a_2 \leq y+A_2$. Hence,
\beao
\dfrac{\PP[X>x-A_1\,, \;Y>y-A_2]}{\PP[X>x\,, \;Y>y]}&\geq& \dfrac{\PP[X>x+a_1\,, \;Y>y+a_2]}{\PP[X>x\,, \;Y>y]}\\[2mm]
&\geq& \dfrac{\PP[X>x+A_1\,, \;Y>y+A_2]}{\PP[X>x\,, \;Y>y]}\,,
\eeao
where the first fraction tends to unity, as $x\wedge y \to \infty$, by Definition \ref{def.KP.1} and the last fraction also tends to unity, as $x\wedge y \to \infty$, after the change of variables $x'=x+A_1$, and $y'=y+A_2$ and by Definition \ref{def.KP.1}.

Definition \ref{def.KP.2.2} provides the insensitivity property in joint distributions, see the uni-variate analogue for example in \cite{foss:korshunov:zachary:2013} or in \cite{Konstantinides:2018}.  

\bde \label{def.KP.2.2}
Let $a_F(x), a_G(y)>0$ for any $x>0,\,y>0$ be two non-decreasing function. We say that the joint distribution ${\bf F}=(F,\,G)$, with right endpoint $r_{{\bf F}}:=(r_F,\,r_G)=(\infty,\,\infty)$, satisfies $(a_F,\,a_G)$-joint insensitivity, if 
\beao
&&\sup_{|a_1|\leq a_F(x),\,|a_2| \leq a_G(y)}\left| \PP[X>x-a_1, Y >y-a_2]-\PP[X>x, Y>y]\right|\\[2mm]
&&\qquad \qquad \qquad =o\left( \PP[X>x,Y>y]\right)\,,
\eeao
as $x\wedge y \to \infty $.
\ede

Now we show that class $ \mathcal{L}^{(2)}$ satisfies the  $(a_F,\,a_G)$-joint insensitive property.

\ble \label{lem.KP.2.2}
Let assume that $(F,\,G) \in \mathcal{L}^{(2)}$. Then there exists some functions $a_F(x),\,a_G(y)$ such that $a_F(x) \to \infty$ and $a_G(y) \to \infty$, as $x\wedge y \to \infty$, and  $(F,\,G)$ satisfies the $(a_F,\,a_G)$-joint insensitive property.
\ele

\pr~
For any integer $n \in \bbn$, from relation \eqref{eq.KP.2.8} we can choose an increasing to infinity sequence $\{u_n\}$, such that the inequality 
\beao
\sup_{|a_1|\leq n,\,|a_2| \leq n}\left| \PP[X>x-a_1,\,Y >y-a_2]-\PP[X>x,\,Y>y]\right|\leq \dfrac {\PP[X>x,\,Y>y]}n\,,
\eeao
holds for any $x\geq u_n$ and any $y\geq u_n$. Without loss of generality we consider that the sequence  $\{u_n\}$ increases to infinity. We put $a_F(x)=a_G(y)=n$, for any $(x,\,y)\in (u_n,\,u_{n+1}]^2$. From the fact that $u_n \to \infty$, as $\nto$,  we obtain that $a_F(x) \to \infty$, as $\xto$, and $a_G(y) \to \infty$, as $\yto$.

So, form the construction of $a(\cdot)$ we conclude that
\beao
\sup_{|a_1|\leq a_F(x),\,|a_2| \leq a_G(y)}\left| \PP[X>x-a_1,\,Y >y-a_2]-\PP[X>x,\,Y>y]\right|\leq \dfrac{\PP[X>x,\,Y>y]}{n}\,,
\eeao 
for any $x > u_n$ and any $y> u_n$, which is the required result.
~\halmos

\bre \label{rem.KP.2.1}
From the $(a_F,\,a_G)$-joint insensitivity does not follow necessarily the $a_F$ and $a_G$ are insensitivity functions for the marginal distributions $F,\,G$ respectively. Furthermore Lemma \ref{lem.KP.2.2} asserts that
\beao
\lim_{x\wedge y \to \infty }\dfrac{ \PP[X>x\pm a_F(x)\,,\;Y >y \pm a_G(y)]}{\PP[X>x\,,\;Y>y]}=1\,.
\eeao
\ere

Let see now two examples, that help either to understanding or to constructing of such bi-variate distributions. In first case is the simplest, as we construct $(X,\,Y) \in \mathcal{L}^{(2)}$ through the independence between $X$ and $Y$.

\bexam \label{exam.KP.1}
Let $X$ and $Y$ be random variables with distributions $F \in \mathcal{L}$ and $G \in \mathcal{L}$ respectively. We assume that $X$ and $Y$ are independent, to obtain
\beao
\lim_{x\wedge y \to \infty} \dfrac{{\bf \bF_1}(x-a_1,\,y-a_2)}{{\bf \bF_1}(x,\,y)}&=&\lim_{x\wedge y \to \infty} \dfrac{\PP[X>x-a_1,\,Y>y-a_2]}{\PP[X>x,\,Y>y]}\\[2mm]
&=&\lim_{x\wedge y \to \infty}\dfrac{\PP[X>x-a_1]}{\PP[X>x]}\,\dfrac{\PP[Y>y-a_2]}{\PP[Y>y]}=1\,.
\eeao
Therefore $(F,\,G) \in \mathcal{L}^{(2)}$.
\eexam

The next example makes sense, as it can not be reduced into uni-variate distributions. The following dependence structure can be found in \cite{li:2018b}. We say that the random variables $X$ and $Y$ are strongly asymptotic independent (SAI) if hold $\PP[X^- > x,\,Y>y]=O[F(-x)\,\bG(y)]$, $\PP[X > x,\,Y^->y]=O[\bF(x)\,G(-y)]$, as $x\wedge y \to \infty$, and there exists a constant $C>0$ such that holds
\beam \label{eq.KP.12}
\PP[X>x,\,Y>y] \sim C\,\bF(x)\,\bG(y)\,,
\eeam
as $x\wedge y \to \infty$.

If the $X$ and $Y$ are bounded from below, then \eqref{eq.KP.12} is enough to be SAI.

\bexam \label{exam.KP.2}
Let $X$ and $Y$ be random variables with strongly asymptotic independence, with some constant $C>0$ and distributions  $F \in \mathcal{L}$ and $G \in \mathcal{L}$ respectively. Then
\beao
\lim_{x\wedge y \to \infty} \dfrac{{\bf \bF_1}(x-a_1,\,y-a_2)}{{\bf \bF_1}(x,\,y)}&=&\lim_{x\wedge y \to \infty} \dfrac{\PP[X>x-a_1,\,Y>y-a_2]}{\PP[X>x,\,Y>y]}\\[2mm]
&=&\lim_{x\wedge y \to \infty}\dfrac{C\,\bF(x-a_1)\,\bG(y-a_2)}{C\,\bF(x)\,\bG(y)}=1\,.
\eeao
Therefore $(F,\,G) \in \mathcal{L}^{(2)} $.
\eexam

The first two examples restrict themselves either in independent case or in some kind of asymptotic independence. Let notice that in the next example as class $ \mathcal{L}^{(2)}$ we understand the class from Definition \ref{def.KP.1}, but restricted, with respect to convergence, instead of $x \wedge y$ to $x=y$ only. In \cite{li:yang:2015} was used the dependence structure from relation \eqref{eq.KP.2.1*}, through the survival copula $\hat{C}$, to depict the dependence relation among claims in a bi-variate, continuous time risk model. 
We assume that for two random variables $X,\,Y$ follow a survival copula $\hat{C}$, there exists some constant $\gamma \geq 1$ and a positive measurable function $h(\cdot,\,\cdot)$, such that the asymptotic relation
\beam \label{eq.KP.2.1*}
\hat{C}(t_1\,x,\,t_2\,x) \sim x^{\gamma}\,h(t_1,\,t_2)\,,
\eeam
as $x \downarrow 0$ holds, for any $(t_1,\,t_2) \in (0,\,\infty)$.

\bexam \label{exam.KP.3.3*}
Let the random variables $X,\,Y$ follow a survival copula from relation \eqref{eq.KP.2.1*} and $F,\,G$ be their marginal distributions. Furthermore, we assume that it holds
\beam \label{eq.KP.2.2*}
\lim_{\xto} \dfrac{\bG(x)}{\bF(x)}=c\,,
\eeam
for some positive constant $c>0$ and either $F\in \mathcal{L} $ or $G\in \mathcal{L} $ is true. Finally, we suppose that relation \eqref{eq.KP.2.1*} holds with $\gamma =1$. Then we obtain $F,\,G \in \mathcal{L}$, that follows from the closure property of class $ \mathcal{L}$ with respect to strong equivalence of \eqref{eq.KP.2.2*} see \cite{leipus:siaulys:konstantinides:2023}. From \cite[Prop. 3.1]{li:yang:2015} we have the  random variables $X,\,Y$ to be asymptotic dependent, and further they satisfy 
\beao
\lim_{\xto} \dfrac{\PP[X>x,\,Y>x]}{\PP[X>x]} =h(1,\,c)> 0\,,
\eeao
hence by the last formulas, for any $(a_1,\,a_2)> (0,\,0)$ it holds
\beao
\lim_{\xto} \dfrac{\PP[X>x-a_1,\,Y>x-a_2]}{\PP[X>x,\,Y>x]} =\lim_{\xto} \dfrac{h(1,\,c)\,\PP[X>x-a_1]}{h(1,\,c)\,\PP[X>x]}=1\,,
\eeao
so we find $(F,\,G) \in \mathcal{L}^{(2)}$, in the sense that in Definition \ref{def.KP.1} the convergence is valid with $x=y$.
\eexam

We can find several dependence structures that satisfy the $\mathcal{L}^{(2)}$ condition. However, we choose to pursue theoretical results.

Now we pass to the bi-variate subexponential distribution class $\mathcal{S}^{(2)}$.

\bde \label{def.KP.2}
We say that the random pair $(X,\,Y)$, with marginal distributions $F$ and $G$ respectively, belongs to the class of bi-variate subexponential distributions, symbolically $(F,\,G) \in \mathcal{S}^{(2)}$, if
\begin{enumerate}
\item
$F \in \mathcal{S}$ and $G \in \mathcal{S}$.
\item
$(F,\,G) \in \mathcal{L}^{(2)} $.
\item
It holds
\beam \label{eq.KP.13}
\lim_{x\wedge y \to \infty} \dfrac{\PP[X_1+X_2>x\,,\;Y_1+Y_2>y]}{\PP[X>x\,,\;Y>y]}=2^2\,,
\eeam
where $(X_1,\,Y_1)$ and $(X_2,\,Y_2)$ are independent and identically distributed copies of $(X,\,Y)$. 
\end{enumerate}
\ede

\bre
In case of $d$-variate distribution relation \eqref{eq.KP.13} becomes 
\beao
\lim_{x\wedge y \to \infty} \dfrac{\PP[X_{1,1}+X_{1,2}>x_1\,,\;\ldots,\,X_{d,1}+X_{d,2}>x_d]}{\PP[X_{1,1}>x_1\,,\;\ldots,\,X_{d,1}>x_d]}=2^d\,.
\eeao
\ere

In the following counterexample we show for what reason the $\mathcal{L}^{(2)}$ property does NOT contain always the property \eqref{eq.KP.13} in Definition \ref{def.KP.2}. We should mention that in some papers the $SAI$ property  \eqref{eq.KP.12} can hold with $C=0$ too, see \cite{li:2018a} and \cite{ji:wang:yan:cheng:2023}, for examples in terms of copulas.

\bexam \label{exam.KP.2.2}
Let $(X,\,Y)$ be a random pair with marginal distributions $(F,\,G)$, such that $F \in \mathcal{R}_{-\alpha_1}$, $G \in \mathcal{R}_{-\alpha_2}$, with $\alpha_1,\,\alpha_2 \in (0,\,\infty)$. Let us assume that $(X,\,Y)$ satisfy $SAI$ property with $C=0$. By \cite[Lem. 3.1]{li:2024}, for degenerated weights to unit and $n=m=2$, we obtain relation \eqref{eq.KP.13}. Furthermore, we see that $F,\,G \in \mathcal{S} \subsetneq \mathcal{L}$,  but for any $(a_1,\,a_2)> (0,\,0)$ the ratio
\beao
1 \leq \lim_{x\wedge y \to \infty}  \dfrac{\PP[X>x-a_1,\,Y>y-a_2]}{\PP[X>x,\,Y>y]}\,, 
\eeao
is NOT necessarily bounded from above by  unit. Hence, relation \eqref{eq.KP.13} and the subexponentialy of the marginals do NOT imply the $\mathcal{L}^{(2)}$ property.
\eexam

This last example replies to the following issue. Weather $F,\,G \in \mathcal{L}$ does not imply directly that the joint distribution tail tends to unit. The example of $SAI$ structure with $C=0$, provides the next statement. If $F, \,G \in \mathcal{B}$ does NOT imply directly that $(F,\,G) \in \mathcal{B}^{(2)}$, where $ \mathcal{B} \in \{ \mathcal{C},\, \mathcal{D},\, \mathcal{L}\}$.

Now we come to the bi-variate dominatedly varying distribution class $\mathcal{D}^{(2)}$.

\bde \label{def.KP.3}
We say that the random pair $(X,\,Y)$, with marginal distributions $F$ and $G$ respectively, belongs to the class of bi-variate dominatedly varying  distributions, symbolically $(F,\,G) \in \mathcal{D}^{(2)}$, if
\begin{enumerate}
\item
$F \in \mathcal{D}$ and $G \in \mathcal{D}$.
\item
It holds
\beam \label{eq.KP.14}
\limsup_{x\wedge y \to \infty} \dfrac{{\bf \bF_b}(x,\,y)}{{\bf \bF_1}(x,\,y)}=\limsup_{x\wedge y \to \infty} \dfrac{\PP[X>b_1\,x,\,Y>b_2\,y]}{\PP[X>x,\,Y>y]} < \infty\,,
\eeam
for some, or equivalently for all ${\bf b}=(b_1,\,b_2) \in (0,\,1)^2$, with $b_1$ not necessarily equal to $b_2$.
\end{enumerate}
\ede

It is obvious that \ref{eq.KP.14} is equivalently with:
\beao
\liminf_{(x,y)\rightarrow(\infty,\infty)}\dfrac{{\bf \bF_b}(x,\, y)}{{\bf \bF_1}(x, y)}>0
\eeao
for some, or equivalently for all ${\bf b}=(b_1,\,b_2) \in (1,\,\infty)^2$, with $b_1$ not necessarily equal to $b_2$.

\bre\label{rem.KP.1}
In \cite{konstantinides:passalidis:2024} was introduced the class $\mathcal{D}_{n}$ (for some $n\in\mathbb{N}$) of multivariate dominatedly varying random vectors. It is obvious that in case $n=2$ our approach include this definition. Namely 
\begin{equation*}
\mathcal{D}_{2}\subset\mathcal{D}^{(2)}	
\end{equation*}
\ere

\bde \label{def.KP.4}
We say that the random pair $(X,\,Y)$, with marginal distributions $F$ and $G$ respectively, belongs to the class of bi-variate consistently varying  distributions, symbolically $(F,\,G) \in \mathcal{C}^{(2)}$, if
\begin{enumerate}

\item
$F \in \mathcal{C}$ and $G \in \mathcal{C}$.

\item
It holds
\beao
\lim_{{\bf z} \uparrow {\bf 1}}\limsup_{x\wedge y \to \infty} \dfrac{{\bf \bF_z}(x,\,y)}{{\bf \bF_1}(x,\,y)}=1\,,
\eeao
or equivalently
\beao
\lim_{{\bf z} \downarrow {\bf 1}}\liminf_{x\wedge y \to \infty} \dfrac{{\bf \bF_z}(x,\,y)}{{\bf \bF_1}(x,\,y)}=1\,,
\eeao
where ${\bf z}=(z_1,\,z_2)$, and ${\bf 1}=(1,\,1)$.
\end{enumerate}
\ede

Examples \ref{exam.KP.1} and \ref{exam.KP.2} remain in tact in classes $\mathcal{D}^{(2)}$ and $\mathcal{C}^{(2)}$, hence they keep functioning in class $(\mathcal{D}\cap \mathcal{L})^{(2)}:=\mathcal{D}^{(2)} \cap \mathcal{L}^{(2)}$.

\bth \label{th.KP.2.1}
It holds $\mathcal{C}^{(2)} \subset \mathcal{L}^{(2)}$.
\ethe

\pr~
Let consider that $(F,\,G) \in \mathcal{C}^{(2)}$.  Then, for ${\bf a}=(a_1,\,a_2)=(0,\,0)$ for any distributions $F,\,G$, we obtain
\beam \label{eq.KP.2.1}
1\leq \liminf_{x \wedge y \to \infty} \dfrac{\PP[X>x-a_1,\,Y>y-a_2]}{\PP[X>x,\,Y>y]}\,.
\eeam
Hence, we have to show that the upper bound of the last fraction is equal to unity. We observe that for any small enough $\delta_1,\,\delta_2 >0$ there exists some $x_0>0$, such that $x\,(1-\delta_1) \leq x- a_1$ and $y\,(1-\delta_2) \leq y- a_2$, for any $x \wedge y \geq x_0$. Therefore we find
\beam \label{eq.KP.2.2}
 \limsup_{x\wedge y \to \infty} \dfrac{\PP[X>x-a_1, Y>y-a_2]}{\PP[X>x, Y>y]}\leq \limsup_{x\wedge y \to \infty} \dfrac{\PP[X>x (1-\delta_1),\,Y>y (1-\delta_2)]}{\PP[X>x,\,Y>y]}\to 1,
\eeam
as $(\delta_1,\,\delta_2) \to (0,\,0)$, where in the last step we use the properties of class $\mathcal{C}^{(2)}$ for the pair of distributions $(F,\,G)$. So, by relations \eqref{eq.KP.2.1} and \eqref{eq.KP.2.2} we conclude that $(F,\,G) \in \mathcal{L}^{(2)}$.  
~\halmos

\section{Max-sum equivalence and closure properties with respect to convolution} \label{sec.KP.3}

Now, we present two definitions. In the first one, we define the closure property with respect to convolution in bi-variate distributions. In this case we formulate the main result, showing that the class $\mathcal{D}^{(2)}$ is closed but the class $\mathcal{L}^{(2)}$ is not. The second definition, given at the end of section,   under concrete dependence structures, also presented later, is fulfilled with respect to   classes $(\mathcal{D}\cap \mathcal{L})^{(2)}$ and $\mathcal{C}^{(2)}$.

\bde \label{def.KP.5}
Let $X_1,\,X_2,\,Y_1,\,Y_2$ be random variables, with distributions $F_1$, $F_2$, $G_1 $and $G_2$ respectively. If the following conditions are true
\begin{enumerate}
\item
$F_1 \in \mathcal{B}$, $F_2 \in \mathcal{B}$, $G_1 \in \mathcal{B}$, $G_2 \in \mathcal{B}$ and for any $k,\,l \in \{1,\,2\}$, holds $(F_k,\,G_l) \in \mathcal{B}^{(2)}$,
\item
Holds $(F_{X_1+X_2},\,G_{Y_1+Y_2}) \in \mathcal{B}^{(2)}$,
\end{enumerate}
where $\mathcal{B}^{(2)}$ is some bi-variate class, defined in section \ref{sec.KP.2}, then we say that the class $\mathcal{B}^{(2)}$ is closed with respect to sum. If $X_1$, $X_2$ are independent random variables, and $Y_1$,\,$Y_2$ are also independent, then we say that $\mathcal{B}^{(2)}$ is closed with respect to convolution, symbolically $(F_1*F_1, \,G_1*G_2) \in \mathcal{B}^{(2)}$.
\ede

In last definition, although the check of $F_1*F_2 \in  \mathcal{B}$, $G_1*G_2 \in  \mathcal{B}$ is implies directly by the uni-variate closure properties, still are NOT implied also the check of $(F_k,\,G_l) \in \mathcal{B}^{(2)}$, for any $k,\,l \in \{1,\,2\}$, see for example the case of $SAI$ with $C=0$, and still is NOT implied that the joint tail of $(X_1+X_2,\,Y_1+Y_2)$ has the desired property of $\mathcal{B}^{(2)}$. Hence, we find out that the dependence structures among the components, play crucial role to the closure properties of bi-variate vectors.

Next we see that the class $\mathcal{D}^{(2)}$ is closed with respect to convolution (of arbitrarily dependent random vectors with arbitrarily non-negative dependent components), under the condition that the point $(1)$ in Definition \ref{def.KP.5} is satisfied.

\bth \label{th.KP.1}
Let non-negative random variables $X_1,\,X_2,\,Y_1,\,Y_2$ with distributions $F_1$, $F_2$, $G_1$ and $G_2$, from class $\mathcal{D}$, respectively. We assume that $(X_k,\,Y_l) \in \mathcal{D}^{(2)}$ for any $k,\,l\,\in \{1,\,2\}$, then $(F_{X_1+X_2},\,G_{Y_1+Y_2}) \in \mathcal{D}^{(2)}$.
\ethe

\pr~
At first, for the first condition of $\mathcal{D}^{(2)}$, we obtain $F_{X_1+X_2} \in \mathcal{D}$ and $G_{Y_1+Y_2} \in \mathcal{D}$, because of Proposition \ref{pr.KP.1}.

Taking into consideration that all the distributions have support the interval $[0,\,\infty)$, by the elementary inequalities
\beao
\PP[X_1+X_2>x] \leq \PP[X_1>x/2]+ \PP[X_2 > x/2]\,,\\[2mm]
\PP[X_1+X_2>x] \geq \dfrac 12 \left(\PP[X_1>x]+ \PP[X_2 > x]\right)\,,
\eeao 
we find
\beao
\PP[X_1+X_2>x,\,Y_1+Y_2>y] &\leq& \PP[X_1>x/2,\,Y_1+Y_2>y]+ \PP[X_2>x/2,\,Y_1+Y_2>y]\\[2mm]
&\leq& \PP\left[X_1>\dfrac x2,\,Y_1>\dfrac y2\right]+ \PP\left[X_1>\dfrac x2,\,Y_2>\dfrac y2\right]\\[2mm]
&&+\PP\left[X_2>\dfrac x2,\,Y_1>\dfrac y2\right]+\PP\left[X_2>\dfrac x2,\,Y_2> \dfrac y2\right]\,,
\eeao
hence
\beam \label{eq.KP.16}
\PP[X_1+X_2>x,\,Y_1+Y_2>y]\leq \sum_{k=1}^2 \sum_{l=1}^2 \PP\left[X_k>\dfrac x2,\,Y_l>\dfrac y2\right]\,.
\eeam
From the other side
\beao
&&\PP[X_1+X_2>x,\,Y_1+Y_2>y] \geq \dfrac {\PP[X_1>x,\,Y_1+Y_2>y]+ \PP[X_2>x,\,Y_1+Y_2>y]}2 \\[2mm]
&&\geq \dfrac {\PP[X_1>x,\,Y_1>y]+ \PP[X_1>x,\,Y_2>y]+\PP[X_2>x,\,Y_1>y]+\PP[X_2>x,\,Y_2>y]}4\,,
\eeao
from where we obtain
\beam \label{eq.KP.17}
\PP[X_1+X_2>x,\,Y_1+Y_2>y]\geq  \dfrac 14\,\sum_{k=1}^2 \sum_{l=1}^2 \PP[X_k>x,\,Y_l>y]\,.
\eeam
Therefore by relations \eqref{eq.KP.16} and \eqref{eq.KP.17}, due to  $(X_k,\,Y_l) \in \mathcal{D}^{(2)}$ for any $k,\,l\,\in \{1,\,2\}$, and ${\bf b}=(b_1,\,b_2) \in (0,\,1)^2$ we find
\beao
&&\limsup_{x\wedge y \to \infty}\dfrac{\PP[X_1+X_2>b_1\,x,\,Y_1+Y_2>b_2\,y]}{\PP[X_1+X_2>x,\,Y_1+Y_2>y]}\\[2mm]
&&\qquad \qquad \leq 4\,\limsup_{x\wedge y \to \infty}\dfrac{\sum_{k=1}^2 \sum_{l=1}^2 \PP[X_k>b_1\,x,\,Y_l>b_2\,y]}{\sum_{k=1}^2 \sum_{l=1}^2 \PP[X_k>x,\,Y_l>y]}\\[2mm]
&&\qquad \qquad \leq 4\,\max_{k,\,l\,\in \{1,\,2\}}\left\{ \limsup_{x\wedge y \to \infty}\dfrac{\PP[X_k>b_1\,x,\,Y_l>b_2\,y]}{\PP[X_k>x,\,Y_l>y]}\right\} < \infty\,.
\eeao
So we conclude $(F_{X_1+X_2},\,G_{Y_1+Y_2}) \in \mathcal{D}^{(2)}$.
~\halmos

\bre
Let notice here that the vectors $(X_k,\,Y_l)$ for $k,\,l =1,\,2$, are NOT necessarily under the same dependence structure, for example we can have $(X_1,\,Y_1)$ with independent components and $(X_1,\,Y_2)$ to be $SAI$, with $C>0$. A case where we see that $(X_k,\,Y_l) \in \mathcal{D}^{(2)}$ for any $k,\,l =1,\,2$ is the following. Let $X_1,\,X_2,\,Y_1,\,Y_2$ with distributions from class $\mathcal{D}$ and also the $X_1,\,X_2$ are arbitrarily dependent and the  $Y_1,\,Y_2$ are arbitrarily dependent. If $(X_k,\,Y_l)$ are $SAI$ with $C_{k,l}>0$, not necessarily the same for each pair, then $(X_k,\,Y_l)\in \mathcal{D}^{(2)}$  for any $k,\,l =1,\,2$.
\ere

Now we are ready to define the max-sum equivalence in two dimensions.

\bde \label{def.KP.6}
Let $X_1,\,X_2,\,Y_1,\,Y_2$ be random variables. Then we say that they are joint max-sum equivalent if
\beao
\PP[X_1+X_2>x,\,Y_1+Y_2>y]\sim \sum_{k=1}^2 \sum_{l=1}^2 \PP[X_k>x,\,Y_l>y]\,,
\eeao
as $x\wedge y \to \infty$.
\ede

This kind of asymptotic relation will be established for classes $(\mathcal{D}\cap \mathcal{L})^{(2)}$ and $\mathcal{C}^{(2)}$, under the assumption of non-negative support and some specific dependence structure.

\section{Joint behavior of random sums} \label{sec.KP.4}

In one dimension, the following asymptotic relation attracted attention
\beam \label{eq.KP.19}
\PP\left[ \sum_{i=1}^n X_i>x \right]\sim \sum_{i=1}^n \PP[X_i>x]\,,
\eeam
as $x\to \infty$. Therefore, we study the behavior of both, the maximum $\bigvee_{i=1}^n X_i$
and the maximum of sums 
\beao
\bigvee_{i=1}^n S_i:=\max_{1\leq  k \leq n}\sum_{i=1}^k X_i\,,
\eeao 
for some distributions and correspondingly with some dependence structures to examine if it holds
\beam \label{eq.KP.20}
\PP\left[\sum_{i=1}^n X_i>x\right] \sim \PP\left[ \bigvee_{i=1}^n X_i>x \right]\sim \PP\left[ \bigvee_{i=1}^n S_i>x \right]\sim \sum_{i=1}^n\PP\left[  X_i>x \right]\,,
\eeam
as $x\to \infty$. The relations \eqref{eq.KP.19} and \eqref{eq.KP.20} have been studied extensively, see for example in \cite{asmussen:foss:korshunov:2003}, \cite{geluk:ng:2006},  \cite{geluk:tang:2009}, \cite{ng:tang:yang:2002}, \cite{jiang:gao:wang:2014}. A similar interest has been appeared for weighted sums of the form
\beao
S_n^{\Theta} := \sum_{i=1}^n \Theta_i\,X_i\,,\qquad \bigvee_{i=1}^n S_i^{\Theta}:=\max_{1\leq k \leq n} \sum_{i=1}^k \Theta_i\,X_i\,,
\eeao 
and for the circumstances when they satisfy relations \eqref{eq.KP.19} and \eqref{eq.KP.20}, see for example \cite{tang:yuan:2014}, \cite{tang:tsitsiashvili:2003b},  \cite{yang:leipus:siaulys:2012}, \cite{zhang:shen:weng:2009}.

In this section we study the relation \eqref{eq.KP.20} in two dimensions. This can be achieved for the class  $(\mathcal{D}\cap \mathcal{L})^{(2)}$ under generalized tail asymptotic dependence. Although the uni-variate randomly weighted sums are well studied, this is not true for the multivariate case.

Let mention some papers, involved in the asymptotic behavior of the joint-tail probability 
\beao
\PP\left[ \sum_{i=1}^n \Theta_i\,X_i >x,\;\sum_{j=1}^n \Delta_i\,Y_i >y \right]\,,
\eeao 
as for example  \cite{chen:yang:2019}, \cite{li:2018b}, \cite{shen:du:2023}, \cite{shen:ge:fu:2018}, \cite{yang:chen:yuen:2024}.
 
We restrict ourselves at moment, in the study of non-weighted  random sums of the following form 
\beao
\PP\left[ \sum_{i=1}^n X_i>x,\;\sum_{j=1}^n Y_j >y \right]\,.
\eeao

We note that, in almost all existing papers, the dependence structure for the main variables $X_i$, $Y_j$ is either of the form: $\{(X_i,\,Y_i)\,,\;i \in \bbn\}$ independent random vectors and there exists some dependence structure in each random pair, or there exists dependence among $X_1,\,\ldots,\,X_n$ and $Y_1,\,\ldots,\,Y_n$, but the $X_i$ and $Y_j$ are independent for any $i,\,j$. Using generalized tail asymptotic independence (GTAI), introduced in \cite{konstantinides:passalidis:2023}, both dependence structures are simultaneously permitted. GTAI is defined as follows. Let consider two sequences of random variables $\{X_n,\;n\in \bbn\}\,,\; \{Y_m,\;m\in \bbn\}$. We say that the random variables $X_1,\,\ldots,\,X_n,\,Y_1,\,\ldots,\,Y_m$ follow the generalized tail asymptotic independence, if 
\begin{enumerate}
\item
It holds 
\beao
\lim_{\min\{x_i,\,x_k,\,y_j\} \to \infty}  \PP[|X_i|> x_i\;|\;X_k>x_k,\, Y_j > y_j]=0\,,
\eeao
for any $1\leq i \neq k \leq n$, $j=1,\,\ldots,\,m$.
\item
It holds 
\beao
\lim_{\min\{x_i,\,y_k,\,y_j\} \to \infty}  \PP[|Y_j|> y_j\;|\;X_i>x_i,\, Y_k > y_k]=0\,,
\eeao
for any $1\leq j \neq k \leq m$, $i=1,\,\ldots,\,n$.
\end{enumerate}
The aim of this dependence structure is the modeling the dependence both, in each sequence of random variables and in the interdependence between the sequences. We have to notice that if the $X_i$ and $Y_j$ are independent for any $i,\,j$, then each sequence of random variables follows tail asymptotic dependence (TAI), see definition bellow, however in any other case the GTAI does not restrict each sequence to TAI, but in a more general form of dependence. 

It is easy to find that GTAI contains the case when $X_1,\,\ldots,\,X_n$ are independent or when $Y_1,\,\ldots,\,Y_n$ are independent or both. Even more this dependence structure indicates that the probability to happen three extreme events, is negligible with respect to probability to happen two extreme events, one in each sequence, and in some sense $GTAI$ belongs to the dependencies of second order asymptotic independence.

In the most of our results we use the TAI dependence structure, as an extra assumption, which characterizes the dependence of the terms of each sequence. This dependence structure was introduced by \cite{geluk:tang:2009}. We say that they are Tail Asymptotic Independent, symbolically $TAI$, and in some works named as strong quasi-asymptotically independent, if for any  pair $i,\,j =1,\,\ldots,\,n$, with $i \neq j$ it holds the limit
\beao 
\lim_{x_i \wedge x_j \to \infty}\PP[ |X_i| >x_i\;|\; X_j >x_j] = 0\,.
\eeao

Next result provides an asymptotic relation for the maximum of two sequences of random variables under the GTAI, WITHOUT to impose some assumption  for the distributions of $
X_1,\,\ldots,\,X_n,\,Y_1,\,\ldots,\,Y_m$, (except the infinite right point).

\bth \label{th.KP.2}
If $X_1,\,\ldots,\,X_n$ are random variables with distributions $F_1,\,\ldots,\,F_n$ respectively, and $Y_1,\,\ldots,\,Y_m$  are random variables with distributions $G_1,\,\ldots,\,G_m$ and $X_1,\,\ldots,\,X_n$, $Y_1,\,\ldots,\,Y_m$ are GTAI then holds
\beao
\PP\left[\bigvee_{i=1}^n X_i>x\,,\;\bigvee_{j=1}^m Y_j>y \right]\sim  \sum_{i=1}^n \sum_{j=1}^m \PP\left[  X_i>x\,,\; Y_j>y \right]\,,
\eeao
as $x\wedge y \to \infty$.
\ethe

\pr~
For $x>0,\,y>0$ holds
\beam \label{eq.KP.23}
\PP\left[\bigvee_{i=1}^n X_i>x\,,\;\bigvee_{j=1}^m Y_j>y \right]\leq  \sum_{i=1}^n \sum_{j=1}^m \PP\left[  X_i>x\,,\; Y_j>y \right]\,.
\eeam
Further for the lower bound we use Bonferroni's inequality 
\beao
&&\PP\left[\bigvee_{i=1}^n X_i>x\,,\;\bigvee_{j=1}^m Y_j>y \right] \\[2mm]
&&\geq  \sum_{i=1}^n \PP\left[  X_i>x\,,\; \bigvee_{j=1}^m Y_j>y \right]- {\sum\sum}_{ i< l =1} ^{n} \PP\left[  X_i>x\,,\;X_l>x\,,\: \bigvee_{j=1}^m Y_j>y \right]  \\[2mm]
&&\geq \sum_{i=1}^n \sum_{j=1}^m  \PP\left[  X_i>x\,,\;  Y_j>y \right]\\[2mm]
&&- \sum_{i=1}^n {\sum\sum}_{j< k =1}^m \PP\left[  X_i>x\,,\;Y_j>y\,,\; Y_k> y \right]\\[2mm]
&&-{\sum\sum}_{l< i=1}^n \sum_{j=1}^m  \PP\left[  X_i>x\,,\;X_l>x\,,\; Y_j>y \right]\\[2mm]
&&=: I_1(x,\,y) - I_2(x,\,y) - I_3(x,\,y)\,.
\eeao
For $I_2(x,\,y)$ we obtain
\beao
I_2(x,\,y)&=&\sum_{i=1}^n {\sum\sum}_{j<k =1}^m \PP\left[  X_i>x\,,\;Y_j>y\,,\; Y_k> y \right]\\[2mm]
&=&\sum_{i=1}^n {\sum\sum}_{j <k=1}^m \PP\left[ Y_k> y\;|\; X_i>x\,,\;Y_j>y \right]\, \PP\left[  X_i>x\,,\;Y_j>y\right]\\[2mm]
&=&o\left( \sum_{i=1}^n \sum_{j=1}^m \PP\left[  X_i>x\,,\;Y_j>y \right]\right)=o[I_1(x,\,y)]\,,
\eeao
as $x\wedge y \to \infty$, where in the last step we use the GTAI property. In a similar way we can find 
\beao
I_3(x,\,y)=o[I_1(x,\,y)]\,,
\eeao 
as $x\wedge y \to \infty$. Hence we conclude
\beam \label{eq.KP.24}
\PP\left[\bigvee_{i=1}^n X_i>x\,,\;\bigvee_{j=1}^m Y_j>y \right]\gtrsim  \sum_{i=1}^n \sum_{j=1}^m \PP\left[  X_i>x\,,\; Y_j>y \right]\,,
\eeam
as $x\wedge y \to \infty$. Now, from relations \eqref{eq.KP.23} and \eqref{eq.KP.24} we have the result.
~\halmos

Before next theorem, we need some preliminary lemmas. Next lemma provides an important property of the GTAI structure, presenting itself as closure property  with respect to sum. 

\ble \label{lem.KP.1}
If $X_1,\,\ldots,\,X_n,\,Y_1,\,\ldots,\,Y_m$ follow the generalized tail asymptotic independence (GTAI), then holds
\beam \label{eq.KP.25}
\lim_{\min\{x_I,\,x_k,\,y_j\} \to \infty} \PP\left[\left| \sum_{i \in I} X_i \right|>x_I\;\Big|\;X_k>x_k\,,\; Y_j>y_j \right] =0\,,
\eeam
for $I\subsetneq \{1,\,\ldots,\,n\}$ and $k \in \{1,\,\ldots,\,n\} \setminus I$, $j=1,\,\ldots,\,m$. Similarly holds
\beam \label{eq.KP.26}
\lim_{\min\{x_i,\,y_k,\,y_J\} \to \infty} \PP\left[\left| \sum_{j \in J} Y_j \right|>y_J\;\Big|\;Y_k>y_k\,,\; X_i>x_i \right] =0\,,
\eeam
for $J\subsetneq \{1,\,\ldots,\,m\}$ and $k \in \{1,\,\ldots,\,m\} \setminus J$, $i=1,\,\ldots,\,n$.
\ele

\pr~
It is enough to show relation \eqref{eq.KP.25} as relation \eqref{eq.KP.26} follows by similar way. Indeed, we observe that
\beao
&&\lim_{\min\{x_I,\,x_k,\,y_j\} \to \infty} \PP\left[\left| \sum_{i \in I} X_i \right|>x_I\;\Big|\;X_k>x_k\,,\; Y_j>y_j \right] \\[2mm]
&&\qquad \qquad \leq \lim_{\min\{x_I,\,x_k,\,y_j\} \to \infty} \sum_{i \in I}\PP\left[\left|  X_i \right|>\dfrac{x_I}n\;\Big|\;X_k>x_k\,,\; Y_j>y_j \right]=0\,,
\eeao
where the last step follows from $GTAI$ property.
~\halmos

In most of following results, we assume that the random variables $X_1,\,\ldots,\,X_n$, $Y_1,\,\ldots,\,Y_n$ are  $GTAI$ and follow distributions from some class $ \mathcal{B} \in \{ \mathcal{C},\, \mathcal{D}\cap \mathcal{L},\, \mathcal{L}\}$, and the same time it holds $(X_k,\,Y_l) \in \mathcal{B}^{(2)}$, for any $k,\,l \in \{1,\,\ldots,\,n\}$. Following a referee's advice, we provide some examples, to show that $GTAI$ and $(X_k,\,Y_l) \in \mathcal{B}^{(2)}$ are possible simultaneously. For sake of simplicity we consider only the case $n=2$ with non-negative random variables.   

\bexam \label{exam.KP.4.1}
Let $X_1,\,X_2,\,Y_1,\,Y_2$ be non-negative random variables with distributions from class $ \mathcal{B} \in \{ \mathcal{C},\, \mathcal{D}\cap \mathcal{L},\, \mathcal{L}\}$. Further, we suppose that the $X_1,\,X_2$ are $TAI$, the $Y_1,\,Y_2$ are also $TAI$, while the $(X_1,\,X_2)$ and $(Y_1,\,Y_2)$ are independent random pairs. Then, we directly find that $(X_k,\,Y_l)  \in \mathcal{B}^{(2)}$ and additionally by $TAI$ we obtain that for any $\vep >0$, there exists some $x_0>0$, such that for any $1\leq i \neq k \leq 2$ it holds $\PP[X_i>x_i\;|\;X_k > x_k] < \varepsilon$, for any $x_i\wedge x_k \geq x_0$. Hence for any $1\leq i \neq k \leq 2$, $j=1,\,2$ we conclude
\beao
\PP\left[X_i>x_i\,,\; X_k> x_k \,,\;Y_j>y_j\right]&=& \PP\left[X_i>x_i\,,\; X_k> x_k\right]\,\PP\left[Y_j>y_j\right] \\[2mm]
&<& \vep\,\PP\left[X_k> x_k\right]\,\PP\left[Y_j>y_j\right] = \vep\,\PP\left[X_k> x_k,\,Y_j>y_j\right]\,,
\eeao
for any  $x_i\wedge x_k \geq x_0$. From the last relation, because of the arbitrarily choice of $\vep$, we get $\PP\left[X_i>x_i\;|\; X_k> x_k \,,\;Y_j>y_j\right] \to 0$, as $x_i\wedge x_k \wedge y_j \to \infty$.

Similarly, by symmetry, we obtain for any $1\leq j \neq k \leq 2$, $i=1,\,2$ the convergence  $\PP\left[Y_j>y_j\;|\; X_i> x_i \,,\;Y_k>y_k\right] \to 0$, as $x_i\wedge y_k \wedge y_j \to \infty$. Hence, the $X_1,\,X_2,\,Y_1,\,Y_2$ satisfy the $GTAI$.
\eexam

\bexam \label{exam.KP.4.2}
Let $X_1,\,X_2,\,Y_1,\,Y_2$ be non-negative random variables with distributions from class $ \mathcal{B} \in \{ \mathcal{C},\, \mathcal{D}\cap \mathcal{L},\, \mathcal{L}\}$. We suppose that $Z_i,\,Z_j,\,Z_k \in \{X_1,\,X_2,\,Y_1,\,Y_2\}$ with $Z_i \neq Z_j \neq Z_k$ and  $z_i,\,z_j,\,z_k \in \{x_1,\,x_2,\,y_1,\,y_2\}$, where the $z_i,\,z_j,\,z_k$ are aloud to be equal. Let assume the $SAI$ property for any duo or trio of them, namely $\PP[Z_i>z_i\,,\;Z_j>z_j] \sim C_{ij}\,\PP[Z_i>z_i]\,\PP[Z_j>z_j]$, as $z_i\wedge z_j \to \infty$, with $C_{ij} >0$ and  $\PP[Z_i>z_i\,,\;Z_j>z_j\,,\;Z_k>z_k] \sim C_{ijk}\,\PP[Z_i>z_i]\,\PP[Z_j>z_j]\,\PP[Z_k>z_k]$, as $z_i\wedge z_j\wedge z_k \to \infty$, with $C_{ijk} >0$. Then, by $SAI$ in duo mode, we obtain $(X_k,\,Y_l)  \in \mathcal{B}^{(2)}$, for $k,\,l \in \{1,\,2\}$ (see for Example \ref{exam.KP.2}, for the case of class $\mathcal{L}$). Next,  by $SAI$ in trio mode, we obtain directly the $GTAI$.
\eexam

From here on, we study only the case $n=m$. In next lemma, we find the lower asymptotic bound of the joint tail of the random sums 
\beao
S_n := \sum_{k=1}^n X_k\,,\qquad \qquad T_n := \sum_{l=1}^n Y_l\,,
\eeao 
when the summands follow distributions with long tails and the $\mathcal{L}^{(2)}$ property is true for any pair of the summands distribution. A similar result, for the uni-dimensional case, can be found in \cite{geluk:tang:2009}, where the dependence structure is TAI (tail asymptotic independence). In the next result we find generalization to two dimensions and furthermore the GTAI assumption. Next, we introduce the notations 
\beao
S_{n,k}:=S_n- X_k\,,\qquad T_{n,l}:=T_n- Y_l\,,
\eeao 
for some $k\in\{1,\,\ldots,\,n\}$ and some $l\in\{1,\,\ldots,\,n\}$. In what follows, we can choose 
\beam \label{eq.KP.27}
a=(a_F,\,a_G):=\left(\min_{1\leq k \leq n}a_{F_k}\,,\;\min_{1\leq l \leq n} a_{G_l} \right)\,, 
\eeam
namely the minimum of all the  joint insensitivity functions, that means that the function $a(\cdot)$ is insensitive for all the distribution pairs $(F_k,\,G_l)$, for $k,\,l \in \{1,\,\ldots,\,n\}$. In what follows, for sake of simplicity, the function $a(\cdot)$ is understood either as $a_{F}$ for the $X_k$ or as $a_{G}$ for the $Y_l$.

\ble \label{lem.KP.2}
Let  $X_1,\,\ldots,\,X_n,\,Y_1,\,\ldots,\,Y_n$ be random variables with distributions  $F_1,\,\ldots,\,F_n$, $G_1,\,\ldots,\,G_n$ from class $\mathcal{L}$ respectively. We also assume that $X_1,\,\ldots,\,X_n,\,Y_1,\,\ldots,\,Y_n$ satisfy the GTAI property and holds 
\beao
(X_k,\,Y_l)\in  \mathcal{L}^{(2)}\,,
\eeao 
for any $k,\,l \in \{1,\,\ldots,\,n\}$. Then holds
\beao
 \PP\left[S_n>x\,,\;T_n>y \right] \gtrsim  \sum_{k=1}^n \sum_{l=1}^n \PP\left[  X_k>x\,,\; Y_l>y \right]\,,
\eeao
as $x\wedge y \to \infty$.
\ele

\pr~
We choose as $a(\cdot)$ a function with joint insensitivity property for any random pair $(X_k,\,Y_l)$ for   any $k,\,l \in \{1,\,\ldots,\,n\}$. A possible choice of this function is by \eqref{eq.KP.27}. Next, we apply twice Bonferroni's inequality to obtain
\beam \label{eq.KP.4.1}
&&\PP\left[S_n> x\,,\;T_n>y \right] \geq \PP\left[S_n> x\,,\;T_n>y\,,\;\bigvee_{k=1}^n X_k> x+a(x)\,,\; \bigvee_{l=1}^n Y_l>y+a(y) \right] \notag \\[2mm]
&&\geq \sum_{k=1}^n \sum_{l=1}^n \PP\left[S_n> x\,,\;T_n>y\,,\; X_k> x+a(x)\,,\; Y_l>y+a(y) \right]  \\[2mm]
&&- \sum \sum_{1\leq k < i \leq n} \sum_{l=1}^n \PP\left[X_i>x+a(x)\,,\; X_k> x+a(x)\,,\; Y_l>y+a(y) \right] \notag \\[2mm] \notag
&&- \sum_{k=1}^n \sum \sum_{1\leq l < i \leq n} \PP\left[ X_k> x+a(x)\,,\; Y_l>y+a(y)\,,\;Y_i>y+a(y) \right] \\[2mm] \notag
&&=:\sum_{i=1}^3 J_i(x,\,y)\,.
\eeam
Now for each term of $J_2(x,\,y)$ we find
\beao
&& \PP\left[X_i>x+a(x)\,,\; X_k> x+a(x)\,,\; Y_l>y+a(y) \right]  \\[2mm] 
&&=\PP\left[ X_i> x+a(x)\;|\; X_k> x+a(x)\,,\;Y_l>y+a(y) \right] \,\PP\left[ X_k> x+a(x)\,,\; Y_l>y+a(y) \right] \\[2mm] 
&&=o(\PP\left[ X_k> x\,,\;Y_l>y \right])\,,
\eeao
as $x\wedge y \to \infty$, that follows from GTAI property, $\mathcal{L}^{(2)}$ membership and the definition of function $a(\cdot)$. So
\beam \label{eq.KP.4.2}
J_2(x,\,y) =o(\PP\left[ X_k> x\,,\;Y_l>y \right])\,,
\eeam
as $x\wedge y \to \infty$. Similarly, due to symmetry, we have
\beam \label{eq.KP.4.3}
J_3(x,\,y) =o(\PP\left[ X_k> x\,,\;Y_l>y \right])\,,
\eeam
as $x\wedge y \to \infty$.

Finally, for the first term we obtain
\beam \label{eq.KP.4.12*}
J_1(x,\,y) &\geq& \sum_{k=1}^n \sum_{l=1}^n \PP[X_k>x+a(x)\,,\;Y_l > y +a(y)]\\[2mm] \notag
&-& \sum_{1\leq k < i \leq n}^n \sum \sum_{l=1}^n \PP\left[X_k>x+a(x)\,,\;Y_l > y +a(y)\,,\; X_i< -\dfrac {a(x)}n\right]\\[2mm] \notag
&-& \sum_{k=1}^n \sum_{1\leq l < i \leq n}^n \sum  \PP\left[X_k>x+a(x)\,,\;Y_l > y +a(y)\,,\; Y_i< -\dfrac {a(y)}n\right]\,,
\eeam
hence, the last two terms in \eqref{eq.KP.4.12*}, from the $GTAI$ structure and the definition of the function $a$, in combination with properties of class $\mathcal{L}^{(2)}$, become negligible with respect to the fist term in \eqref{eq.KP.4.12*}. Therefore it holds
\beam \label{eq.KP.4.13*}
J_1(x,\,y) &\gtrsim& \sum_{k=1}^n \sum_{l=1}^n \PP[X_k>x+a(x)\,,\;Y_l > y +a(y)]\,,
\eeam
as $x\wedge y \to \infty$.Thus, form relations \eqref{eq.KP.4.2}, \eqref{eq.KP.4.3} and \eqref{eq.KP.4.13*}, together with relation \eqref{eq.KP.4.1} render the desired lower bound. 
~\halmos

\ble \label{lem.KP.4.3*}
Let $X_1,\,X_2,\,Y_1,\,Y_2$ be non-negative random variables, with $GTAI$ property, such that the pair $(X_k,\,Y_l) \in (\mathcal{D}\cap \mathcal{L})^{(2)}$, for any $k,\,l \in \{1,\,2\}$. If by $a$ we denote the insensitivity function  from \eqref{eq.KP.27}, then it holds
\beam \label{eq.KP.4.1*} \notag
 \PP\left[X_1>x-a(x)\,,\;Y_1 >a(y)\,,\; Y_2\leq \dfrac {y}2\right] &\sim& \PP\left[X_1>x\,,\;Y_1 >a(y)\,,\; Y_2\leq \dfrac {y}2\right]\\[2mm]
&\sim& \PP[X_1>x\,,\;Y_1 >a(y)]\,,
\eeam
as $x\wedge y \to \infty$ and further it holds
\beao
 \PP\left[Y_1>y-a(y)\,,\;X_1 >a(x)\,,\; X_2\leq \dfrac {x}2\right] &\sim& \PP\left[Y_1>y\,,\;X_1 >a(x)\,,\; X_2\leq \dfrac {x}2\right]\\[2mm]
&\sim& \PP[Y_1>y\,,\;X_1 >a(x)]\,,
\eeao
as $x\wedge y \to \infty$.
\ele

\pr~
We show only the first relation \eqref{eq.KP.4.1*}, since the second follows along similar steps, due to symmetry. At first, by definition of insensitivity function $a$ from \eqref{eq.KP.27}, and if the pair $(X,\,Y) \in (\mathcal{D}\cap \mathcal{L})^{(2)}$ we obtain
\beam \label{eq.KP.4.3*}
1\leq \lim_{x\wedge y \to \infty} \dfrac{\PP\left[X>x-a(x),\;Y >y\right]}{\PP\left[X>x\,,\;Y >y\right]} \leq \lim_{x\wedge y \to \infty} \dfrac{\PP\left[X>x-a(x),\;Y >y-a(y)\right]}{\PP\left[X>x,\;Y >y\right]}=1,
\eeam
Hence, because of $(X_1,\,Y_1) \in (\mathcal{D}\cap \mathcal{L})^{(2)} \subsetneq \mathcal{L}^{(2)}$, we have through \eqref{eq.KP.4.3*} that it holds
\beao
\PP\left[X_1>x-a(x),\;Y_1 > a(y) \right]\sim \PP\left[X_1>x\,,\;Y_1 >a(y)\right]\,,
\eeao
or equivalently
\beam \label{eq.KP.4.4*} \notag
&&\PP\left[X_1>x-a(x),\;Y_1 > a(y),\;Y_2 \leq \dfrac y2 \right] +\PP\left[X_1>x-a(x),\;Y_1 > a(y),\;Y_2 > \dfrac y2\right]\\[2mm]
&&\sim \PP\left[X_1>x,\;Y_1 > a(y),\;Y_2 \leq \dfrac y2 \right] +\PP\left[X_1>x,\;Y_1 > a(y),\;Y_2 > \dfrac y2\right]\,,
\eeam
as $x\wedge y \to \infty$. We compare the first terms of each side of \eqref{eq.KP.4.4*}, to find
\beao
&&\PP\left[X_1>x-a(x),\;Y_1 > a(y),\;Y_2 \leq \dfrac y2 \right]= \PP\left[X_1>x-a(x),\;Y_1 > a(y)\right]\\[2mm] 
&&-\PP\left[X_1>x-a(x),\;Y_1 > a(y),\;Y_2 > \dfrac y2 \right] =\PP\left[X_1>x-a(x),\;Y_1 > a(y)\right]\\[2mm] 
&&-\PP\left[Y_2>\dfrac y2 \;|\;X_1>x-a(x),\;Y_1 > a(y) \right] \,\PP\left[X_1>x-a(x),\;Y_1 > a(y)\right]\\[2mm] 
&&\sim \PP\left[X_1>x,\;Y_1 > a(y) \right] -o\left(\PP\left[X_1>x,\;Y_1 > a(y)\right]\right)\,,
\eeao
as $x\wedge y \to \infty$, where in the pre-last step we used the class $(\mathcal{D}\cap \mathcal{L})^{(2)}$ property, relation \eqref{eq.KP.4.3*} and the $GTAI$ property. Similarly, we get
\beam \label{eq.KP.4.6*}\notag
&&\PP\left[X_1>x,\;Y_1 > a(y),\;Y_2 \leq \dfrac y2 \right]= \PP\left[X_1>x,\;Y_1 > a(y)\right]\\[2mm] \notag
&&-\PP\left[X_1>x,\;Y_1 > a(y),\,Y_2>  \dfrac y2\right]=\PP\left[X_1>x,\;Y_1 > a(y)\right]\\[2mm] 
&&-\PP\left[Y_2> \dfrac y2 \;\big|\;X_1>x,\;Y_1 > a(y)\right]\PP\left[X_1>x,\;Y_1 > a(y)\right]\\[2mm] \notag
&&\sim \PP\left[X_1>x,\;Y_1 > a(y) \right] -o\left(\PP\left[X_1>x,\;Y_1 > a(y)\right]\right)\,,
\eeam
as $x\wedge y \to \infty$. Therefore, considering all together relations \eqref{eq.KP.4.4*} - \eqref{eq.KP.4.6*}, we conclude relation \eqref{eq.KP.4.1*}.
~\halmos

\bre \label{rem.KP.4.5}
Taking into account relation \eqref{eq.KP.4.4*} - \eqref{eq.KP.4.6*}, together with the fact that the $X_1,\,X_2,\,Y_1,\,Y_2$ are non-negative random variables, which are $GTAI$, it follows that
\beao
\PP\left[X_1>x-a(x),\;Y_1 > a(y),\;Y_2 > \dfrac y2 \right]= o(\PP\left[X_1>x,\;Y_1 > a(y) \right] )\,,\\[2mm]
\PP\left[Y_1>y-a(y),\;X_1 > a(x),\;X_2 > \dfrac x2 \right]= o(\PP\left[X_1> a(x),\;Y_1 > y \right] )\,,
\eeao 
as $x \wedge y \to \infty$.
\ere

The next result shows that in the non-negative  part of class $(\mathcal{D}\cap \mathcal{L})^{(2)}$ the property of joint max-sum equivalence as also under an extra assumption the closure property with respect to convolution are satisfied, as soon as the GTAI holds.

\ble \label{lem.KP.3}
Let $X_1,\,X_2,\,Y_1,\,Y_2$ be non-negative random variables, with the following distributions $F_1,\,F_2,\,G_1,\,G_2$ from class $\mathcal{D}\cap \mathcal{L}$ respectively. Further we assume that the random variables $X_1,\,X_2,\,Y_1,\,Y_2$ satisfy the GTAI and 
\beao
(X_k,\,Y_l)\in  (\mathcal{D}\cap \mathcal{L})^{(2)}\,,
\eeao 
for any $k,\,l \in \{1,\,2\}$ properties. Then it holds 
\beam \label{eq.KP.36}
\PP\left[X_1+X_2> x\,,\;Y_1+Y_2>y \right] \sim  \sum_{k=1}^2  \sum_{l=1}^2 \PP\left[X_k> x\,,\;Y_l>y \right] \,,
\eeam
as $x\wedge y \to \infty$. 
If further $X_1,X_2$ are TAI and $Y_1, Y_2$ are TAI then
\beao
(X_1+X_2,\,Y_1+Y_2) \in (\mathcal{D}\cap \mathcal{L})^{(2)}\,,
\eeao 

\ele

\pr~
From Lemma \ref{lem.KP.2} and the fact that $(\mathcal{D}\cap \mathcal{L})^{(2)} \subsetneq \mathcal{L}^{(2)}$ we find
\beam \label{eq.KP.40}
\PP[X_1+X_2> x\,,\;Y_1+Y_2>y ]\gtrsim \sum_{k=1}^2  \sum_{l=1}^2 \PP\left[X_k> x\,,\;Y_l>y \right]\,,
\eeam
as $x\wedge y \to \infty$, which provides the lower asymptotic bound. 

Let us examine now the upper asymptotic bound
\beam \label{eq.KP.41} \notag
&&\PP[X_1+X_2> x\,,\;Y_1+Y_2>y ]\leq  \PP\left[X_1> x-a(x)\,,\;Y_1+Y_2>y \right]\\[2mm] \notag
&&+\PP\left[X_2> x-a(x)\,,\;Y_1+Y_2>y \right]+\PP\left[X_1> a(x)\,,\;X_2>\dfrac x2\,,\;Y_1+Y_2>y \right]\\[2mm] \notag
&&+\PP\left[X_1> \dfrac x2\,,\;X_2>a(x)\,,\;Y_1+Y_2>y \right] \leq \PP\left[X_1> x-a(x)\,,\;Y_1>y-a(y) \right]\\[2mm] \notag
&&+\PP\left[X_1> x-a(x)\,,\;Y_2>y-a(y) \right]+\PP\left[X_1> x-a(x)\,,\;Y_1>a(y)\,,\;Y_2>\dfrac y2 \right]\\[2mm] \notag
&&+\PP\left[X_1> x-a(x)\,,\;Y_1>\dfrac y2\,,\;Y_2>a(y) \right]+\PP\left[X_2> x-a(x)\,,\;Y_1>y-a(y) \right]\\[2mm] \notag
&&+\PP\left[X_2> x-a(x)\,,\;Y_2>y-a(y) \right]+\PP\left[X_2> x-a(x)\,,\;Y_1>a(y)\,,\;Y_2>\dfrac y2 \right]\\[2mm] \notag
&&+\PP\left[X_2> x-a(x)\,,\;Y_1>\dfrac y2\,,\;Y_2>a(y) \right]+\PP\left[X_1> a(x)\,,\;X_2>\dfrac x2\,,\;Y_1>y-a(y)\right]\\[2mm] \notag
&&+\PP\left[X_1> a(x)\,,\;X_2>\dfrac x2\,,\;Y_2>y-a(y) \right]\\[2mm] \notag
&&+\PP\left[X_1> a(x)\,,\;X_2>\dfrac x2\,,\;Y_1>a(y)\,,\;Y_2>\dfrac y2 \right]\\[2mm] \notag
&&+\PP\left[X_1> a(x)\,,\;X_2>\dfrac x2\,,\;Y_1>\dfrac y2\,,\;Y_2>a(y) \right]\\[2mm] \notag
&&+\PP\left[X_1>\dfrac x2\,,\;X_2> a(x)\,,\;Y_1>y-a(y)\right]+\PP\left[X_1>\dfrac x2\,,\;X_2> a(x)\,,\;Y_2>y-a(y) \right]\\[2mm] \notag
&&+\PP\left[X_1>\dfrac x2\,,\;X_2> a(x)\,,\;Y_1>a(y)\,,\;Y_2>\dfrac y2 \right]\\[2mm]
&&+\PP\left[X_1>\dfrac x2 \,,\;X_2>a(x)\,,\;Y_1>\dfrac y2\,,\;Y_2>a(y) \right]=:\sum_{i=1}^{16}I_i(x,\,y)\,.
\eeam
Taking into account the property $\mathcal{L}^{(2)}$ and the definition of function $a(x)$ we find the asymptotic expressions for $
I_1(x,\,y)\sim \PP\left[X_1> x\,,Y_1>y\right]$, $I_2(x,\,y)\sim \PP\left[X_1> x\,,Y_2>y\right]$, $I_5(x,\,y)\sim \PP\left[X_2> x\,,Y_1>y\right]$,  $I_6(x,\,y)\sim \PP\left[X_2> x\,,Y_2>y\right]$, as $x\wedge y \to \infty$. Hence
\beam \label{eq.KP.42}
I_1(x,\,y)+I_2(x,\,y)+I_5(x,\,y)+I_6(x,\,y) \sim \sum_{k=1}^2  \sum_{l=1}^2 \PP\left[X_k> x\,,\;Y_l>y \right]\,,
\eeam
as $x\wedge y \to \infty$. 

Next, we follow a similar approach for $I_3(x,\,y)$, $I_4(x,\,y)$, $I_7(x,\,y)$, $I_8(x,\,y)$, $I_9(x,\,y)$, $I_{10}(x,\,y)$, $I_{13}(x,\,y)$ and $I_{14}(x,\,y)$. Now, we obtain by Lemma \ref{lem.KP.4.3*} 
\beao
&&I_3(x,\,y) \sim \PP\left[X_1> x\,,\;Y_1> a(y)\,,\;Y_2>\dfrac y2 \right]\\[2mm]
&&\qquad =\PP\left[Y_1> a(y)\;\Big|\;X_1>x\,,\;Y_2>\dfrac y2 \right]\,\PP\left[X_1> x\,,\;Y_2>\dfrac y2 \right]=o\left(\PP\left[X_2> x\,,Y_2>y\right] \right)\,,
\eeao 
as $x\wedge y \to \infty$, that follows because of properties $(\mathcal{D}\cap \mathcal{L})^{(2)}$ and GTAI.

In similar way we find $I_4(x,\,y)=o(\PP\left[X_1> x\,,Y_2>y\right])$,  $I_7(x,\,y)=o(\PP\left[X_2> x\,,Y_2>y\right])$, $I_8(x,\,y)=o(\PP\left[X_2> x\,,Y_2>y\right])$, $I_9(x,\,y)=o(\PP\left[X_2> x\,,Y_1>y\right])$ and finally $I_{10}(x,\,y)=o(\PP\left[X_1> x\,,Y_2>y\right])$, as $x\wedge y \to \infty$. Hence
\beam \label{eq.KP.43}
\sum_{j=3}^4 I_j(x,\,y)+ \sum_{i=7}^{10}I_i(x,\,y) +I_{13}(x,\,y) +I_{14}(x,\,y)  =o\left(\sum_{k=1}^2  \sum_{l=1}^2 \PP\left[X_k> x ,\,Y_l>y \right]\right),
\eeam
as $x\wedge y \to \infty$.

The $I_{11}(x,\,y)$, $I_{12}(x,\,y)$, $I_{15}(x,\,y)$, $I_{16}(x,\,y)$, can be handled also  similarly
\beao
I_{11}(x,\,y) &\leq& \PP\left[X_2> \dfrac x2\,,\;Y_1>a(y)\,,\;Y_2>\dfrac y2 \right]\\[2mm]
&=&\PP\left[Y_1>a(y)\;|\;X_2> \dfrac x2\,,\;Y_2>\dfrac y2 \right]\,\PP\left[X_2> \dfrac x2\,, Y_2>\dfrac y2 \right]\,,
\eeao
or equivalently $I_{11}(x,\,y) =o(\PP\left[X_2> x\,, Y_2>y \right])$, as $x\wedge y \to \infty$, which follows because of properties $(\mathcal{D}\cap \mathcal{L})^{(2)}$ and GTAI. Similarly we find $I_{1j}(x,\,y)=o(\PP\left[X_k>  x\,, Y_l>y \right])$, for some $k,\,l \in \{1,\,2\}$ and for any $j\in \{1,\,2,\,5,\,6\}$. Therefore we obtain
\beam \label{eq.KP.44}
I_{1j}(x,\,y) =o\left(\sum_{k=1}^2  \sum_{l=1}^2 \PP\left[X_k> x\,,\;Y_l>y \right]\right)\,,
\eeam
as $x\wedge y \to \infty$, for any $j\in \{1,\,2,\,5,\,6\}$.

From \eqref{eq.KP.42}, \eqref{eq.KP.43} and \eqref{eq.KP.44}, in combination with \eqref{eq.KP.41} we find that
\beao
\PP[X_1+X_2> x\,,\;Y_1+Y_2>y ] \lesssim \sum_{k=1}^2  \sum_{l=1}^2 \PP\left[X_k> x\,,\;Y_l>y \right]
\eeao
as $x\wedge y \to \infty$, which in combination with \eqref{eq.KP.40} leads to \eqref{eq.KP.36}. 

Now we check the validity of relation $(X_1+X_2,\,Y_1+Y_2) \in (\mathcal{D}\cap \mathcal{L})^{(2)}$. At first, by \eqref{eq.KP.36} we obtain
\beao
&&\limsup_{x\wedge y \to \infty}\dfrac{\PP\left[X_1+X_2> b_1\,x\,,\;Y_1+Y_2>b_2\,y \right]}{\PP\left[X_1+X_2> x\,,\;Y_1+Y_2>y \right]}=\limsup_{x\wedge y \to \infty}\\[2mm]
&&\dfrac{ \sum_{k=1}^2  \sum_{l=1}^2 \PP\left[X_k> b_1\,x\,,\;Y_l>b_2\,y \right]}{ \sum_{k=1}^2  \sum_{l=1}^2 \PP\left[X_k> x\,,\;Y_l>y \right]} \leq \max_{k,\,l \in \{1, 2 \}} \left\{\limsup_{x\wedge y \to \infty}\dfrac{\PP\left[X_k> b_1 x ,\,Y_l>b_2 y \right]}{\PP\left[X_k> x,\,Y_l>y \right]} \right\} < \infty\,,
\eeao
for any ${\bf b}=(b_1,\,b_2) \in (0,\,1)^2$, this means that, we have the one of two conditions of the closure property with respect to $\mathcal{D}^{(2)}$.

Next, we check the closure property with respect to $\mathcal{L}^{(2)}$. From \eqref{eq.KP.36} we obtain
\beao
&&\limsup_{x\wedge y \to \infty}\dfrac{\PP\left[X_1+X_2> x-a_1\,,\;Y_1+Y_2>y-a_2 \right]}{\PP\left[X_1+X_2> x\,,\;Y_1+Y_2>y \right]}\\[2mm]
&&\qquad \qquad =\limsup_{x\wedge y \to \infty}\dfrac{\sum_{k=1}^2  \sum_{l=1}^2 \PP\left[X_k> x-a_1\,,\;Y_l>y-a_2 \right]}{\sum_{k=1}^2  \sum_{l=1}^2 \PP\left[X_k> x\,,\;Y_l>y \right]}\,,
\eeao
for  any ${\bf a}=(a_1,\,a_2) > {\bf 0}$, and therefore
\beao
&&\limsup_{x\wedge y \to \infty}\dfrac{\PP\left[X_1+X_2> x-a_1\,,\;Y_1+Y_2>y-a_2 \right]}{\PP\left[X_1+X_2> x\,,\;Y_1+Y_2>y \right]}\\[2mm]
&&\qquad \qquad \leq \max_{k,\,l \in \{1,\,2 \}}\left\{ \limsup_{x\wedge y \to \infty}\dfrac{ \PP\left[X_k> x-a_1\,,\;Y_l>y-a_2 \right]}{\PP\left[X_k> x\,,\;Y_l>y \right]} \right\}=1\,,
\eeao
and always
\beao
\liminf_{x\wedge y \to \infty}\dfrac{\PP\left[X_1+X_2> x-a_1\,,\;Y_1+Y_2>y-a_2 \right]}{\PP\left[X_1+X_2> x\,,\;Y_1+Y_2>y \right]}\geq 1\,,
\eeao
that means, we have the one of two conditions of the closure property with respect to $\mathcal{L}^{(2)}$ true. So by the extra assumption of TAI between $X_1,X_2$ and $Y_1,Y_2$ by Lemma 4.1 of \cite{geluk:tang:2009} we have that $X_1+X_2\in\mathcal{D}\cap\mathcal{L}$ and $Y_1+Y_2\in\mathcal{D}\cap\mathcal{L}$, as a result we conclude $(X_1+X_2,\,Y_1+Y_2) \in (\mathcal{D}\cap \mathcal{L})^{(2)}$.~\halmos

An easy example,  where we combine the $GTAI$ property for the  $X_1,\,\ldots,\,X_n$, $Y_1,\,\ldots,\,Y_m$ with the $TAI$ property for each sequence, as for the case of  $X_1,\,\ldots,\,X_n$ to be $TAI$ and $Y_1,\,\ldots,\,Y_m$ to be $TAI$, is found in case of each sequence to be $TAI$, but the two sequences to be independent.

Next, we provide a corollary, following from Lemma \ref{lem.KP.3}, where we establish the closure property with respect to $\mathcal{C}^{(2)}$ and the joint max-sum equivalence, under condition GTAI.

\bco \label{cor.KP.2}
Let $X_1,\,X_2,\,Y_1,\,Y_2$ be non-negative random variables, with the distributions $F_1,\,F_2,\,G_1,\,G_2$ from class $\mathcal{C}$ respectively and they satisfy the GTAI condition. If it holds  $(X_k,\,Y_l)\in  \mathcal{C}^{(2)}$, for any $k,\,l \in \{1,\,2\}$, then  
\beam \label{eq.KP.45}
\PP\left[X_1+X_2> x\,,\;Y_1+Y_2>y \right]\sim \sum_{k=1}^2  \sum_{l=1}^2 \PP\left[X_k> x\,,\;Y_l>y \right]\,,
\eeam
as $x\wedge y \to \infty$. If further $X_1,X_2$ are TAI and $Y_1, Y_2$ are TAI then holds $(X_1+X_2\,,\;Y_1+Y_2)\in  \mathcal{C}^{(2)}$.
\eco 

\pr~
Relation \eqref{eq.KP.45} follows from the fact that $\mathcal{C}^{(2)} \subsetneq (\mathcal{D}\cap \mathcal{L})^{(2)}$ and by application of Lemma \ref{lem.KP.3}.

Next, we check the closure property with respect to convolution. From \eqref{eq.KP.45} we obtain
\beao
\PP[X_1+X_2> d_1\,x\,,\;Y_1+Y_2> d_2\,y ] \sim \sum_{k=1}^2  \sum_{l=1}^2 \PP\left[X_k> d_1\,x\,,\;Y_l>d_2\,y \right]\,,
\eeao
as $x\wedge y \to \infty$, for any  ${\bf d}=(d_1,\,d_2) \in (0,\,1)^2$. Hence
\beao
&&\limsup_{x\wedge y \to \infty}\dfrac{\PP\left[X_1+X_2> d_1\,x\,,\;Y_1+Y_2>d_2\,y \right]}{\PP\left[X_1+X_2> x\,,\;Y_1+Y_2>y \right]}\\[2mm]
&&\qquad \qquad =\limsup_{x\wedge y \to \infty}\dfrac{\sum_{k=1}^2  \sum_{l=1}^2 \PP\left[X_k> d_1\,x\,,\;Y_l>d_2\,y \right]}{\sum_{k=1}^2  \sum_{l=1}^2 \PP\left[X_k> x\,,\;Y_l>y \right]} \\[2mm]
&&\qquad \qquad \leq \limsup_{x\wedge y \to \infty}\max_{k,\,l \in \{1,\,2 \}}\left\{ \dfrac{ \PP\left[X_k> d_1\,x\,,\;Y_l>d_2\,y \right]}{\PP\left[X_k> x\,,\;Y_l>y \right]} \right\}\,,
\eeao
Thus, because of the definition of $\mathcal{C}^{(2)}$ we get
\beao
1&\leq& \lim_{{\bf d} \uparrow {\bf 1}} \limsup_{x\wedge y \to \infty}\dfrac{\PP\left[X_1+X_2> d_1\,x\,,\;Y_1+Y_2>d_2\,y \right]}{\PP\left[X_1+X_2> x\,,\;Y_1+Y_2>y \right]}\\[2mm]
&\leq& \lim_{{\bf d} \uparrow {\bf 1}} \limsup_{x\wedge y \to \infty} \max_{k,\,l \in\{1,\,2\}} \left(\dfrac{\PP\left[X_k> d_1\,x\,,\;Y_l>d_2\,y \right]}{\PP\left[X_k> x\,,\;Y_l>y \right]}\right)\\[2mm]
&\leq&  \max_{k,\,l \in\{1,\,2\}} \left(\lim_{{\bf d} \uparrow {\bf 1}} \limsup_{x\wedge y \to \infty}\dfrac{\PP\left[X_k> d_1\,x\,,\;Y_l>d_2\,y \right]}{\PP\left[X_k> x\,,\;Y_l>y \right]}\right)=1\,,
\eeao
this means that the one of two conditions of closedness under convolution holds.
By the assumptions of TAI in each sequence, and by $\mathcal{C}\subsetneq\mathcal{D}\cap\mathcal{L}$, we use Lemma 4.1 of \cite{geluk:tang:2009} and we take 
\beao
1&\leq& \lim_{{ d_1} \uparrow { 1}} \limsup_{x\to \infty}\dfrac{\PP\left[X_1+X_2> d_1\,x \right]}{\PP\left[X_1+X_2> x \right]}=\lim_{{ d_1} \uparrow { 1}} \limsup_{x \to \infty}\dfrac{\PP\left[X_1> d_1\,x \right]+\PP\left[X_1> d_1\,x \right]}{\PP\left[X_1> x \right]+\PP\left[X_1> x \right]}\\[2mm]
&\leq& \lim_{{ d_1} \uparrow { 1}} \limsup_{x \to \infty} \max_{k \in\{1,\,2\}} \left(\dfrac{\PP\left[X_k> d_1\,x \right]}{\PP\left[X_k> x \right]}\right)\leq \max_{k \in\{1,\,2\}} \left(\lim_{{ d_1} \uparrow { 1}} \limsup_{x \to \infty}\dfrac{\PP\left[X_k> d_1\,x \right]}{\PP\left[X_k> x \right]}\right)=1\,,
\eeao
which gives that $(X_1+X_2)\in\mathcal{C}$. With the same argument we have : $(Y_1+Y_2)\in\mathcal{C}$.
That means $(X_1+X_2\,,\;Y_1+Y_2)\in  \mathcal{C}^{(2)}$. 
~\halmos

Now we can give the main result, where we find an analogue to relation \eqref{eq.KP.20} in two dimensions. 

\bth \label{th.KP.3}
Let $X_1,\,\ldots,\,X_n,\,Y_1,\,\ldots,\,Y_n$ be random variables with the following distributions $F_1,\,\ldots,\,F_n\,G_1,\,\ldots,\,G_n$ from class $\mathcal{D}\cap \mathcal{L}$ respectively and they satisfy the GTAI condition, with $(X_k,\,Y_l)\in (\mathcal{D}\cap\mathcal{L})^{(2)}$, for any $k,\,l \in \{1,\,\ldots,\,n\}$. If further $X_1,\,\ldots,\,X_n$ are TAI and $Y_1,\,\ldots,\,Y_n$ are TAI, then
\beao
\PP\left[ \sum_{k=1}^n X_k>x\,,\;\sum_{l=1}^n Y_l >y\right] &\sim& \PP\left[ \bigvee_{i=1}^n S_i>x\,,\;\bigvee_{j=1}^n T_j >y\right] \\[2mm]
&\sim& \PP\left[ \bigvee_{k=1}^n X_k>x\,,\;\bigvee_{l=1}^n Y_l >y\right] \sim \sum_{k=1}^n  \sum_{l=1}^n \PP\left[X_k> x\,,\;Y_l>y \right],
\eeao
as $x\wedge y \to \infty$.
\ethe

\pr~
By Lemma \ref{lem.KP.2} we find
\beao
 \PP\left[\sum_{k=1}^n X_k>x\,,\;\sum_{l=1}^n Y_l>y \right] \gtrsim  \sum_{k=1}^n \sum_{l=1}^n \PP\left[  X_k>x\,,\; Y_l>y \right]\,,
\eeao
as $x\wedge y \to \infty$. Because of closure property of $(\mathcal{D}\cap\mathcal{L})^{(2)}$ with respect to convolution in the positive part, under GTAI condition, we can apply Lemmas  \ref{lem.KP.3}, \ref{lem.KP.1} and employing induction, we find
\beao
\PP\left[\sum_{k=1}^n X_k>x,\,\sum_{l=1}^n Y_l>y \right] \leq \PP\left[\sum_{k=1}^n X_k^+>x,\,\sum_{l=1}^n Y_l^+>y \right]\sim \sum_{k=1}^n \sum_{l=1}^n \PP\left[X_k>x,\,Y_l>y \right],
\eeao 
as $x\wedge y \to \infty$. Now, taking into consideration Theorem \ref{th.KP.2} we find
\beao
 \PP\left[\sum_{k=1}^n X_k>x,\;\sum_{l=1}^n Y_l>y \right] \sim  \sum_{k=1}^n \sum_{l=1}^n \PP\left[  X_k>x,\; Y_l>y \right] \sim  \PP\left[ \bigvee_{k=1}^n X_k>x,\;\bigvee_{l=1}^n Y_l >y\right],
\eeao 
as $x\wedge y \to \infty$. Finally, due to the inequality
\beao
\PP\left[\sum_{k=1}^n X_k>x ,\,\sum_{l=1}^n Y_l>y \right] \leq\PP\left[ \bigvee_{i=1}^n S_i>x ,\,\bigvee_{j=1}^n T_j >y\right]\leq \PP\left[\sum_{k=1}^n X_k^+>x,\,\sum_{l=1}^n Y_l^+>y \right],
\eeao 
we get the asymptotic relation
\beao
 \PP\left[ \bigvee_{i=1}^n S_i>x\,,\;\bigvee_{j=1}^n T_j >y\right] \sim  \sum_{k=1}^n \sum_{l=1}^n \PP\left[  X_k>x\,,\; Y_l>y \right]\,,
\eeao 
as $x\wedge y \to \infty$.
~\halmos

Recently more and more researchers study two dimensional risk models, we refer to the reader \cite{hu:jiang:2013},  \cite{cheng:yu:2019}, \cite{cheng:2021} among  others.
For 
\beao
U_1(k,\,x):=x - \sum_{i=1}^{k}X_i\,, \qquad U_2(k,\,y):=y -\sum_{j=1}^{k}Y_j\,,
\eeao 
for $1\leq k\leq n$, we define now two ruin times,
\beao
T_{max}:=\inf\left\{1\leq k\leq n\;:\; U_1(k,\,x)\wedge U_2(k,\,y) <0 \right\}\,,
\eeao 
that denotes the first moment when both portfolios are found with negative surplus, and for each portfolio we define:
\beao
T_{1}(x):=\inf\left\{1\leq k\leq n: U_{1}(k,\,x)<0 |U_{1}(0,\,x)=x\right\}\,,
\eeao 
\beao
T_{2}(y):=\inf\left\{1\leq k\leq n: U_{2}(k,\,y)<0 |U_{2}(0,\,y)=y\right\}\,,
\eeao
as a result the second type of ruin type is:
\beao
T_{and}:=\max\left\{T_{1}(x),T_{2}(y)\right\}\,,
\eeao 
that corresponds to the first moment, when both portfolios have been with negative surplus, but not necessarily simultaneously. Hence we define the ruin probabilities as
\beam \label{eq.KP.50}
\psi_{max}(x,\,y,\,n) = \PP[T_{max} \leq n]\,,\qquad \psi_{and}(x,\,y,\,n) = \PP[T_{and} \leq n]\,,
\eeam
for any $n \in \bbn$ and $x,\,y >0$. From \eqref{eq.KP.50} we easily find out that
\beao
\psi_{and}(x,\,y,\,n) = \PP\left[ \bigvee_{i=1}^n S_i >x\,,\;\bigvee_{i=1}^n T_i >y\right]\,.
\eeao
Therefore, by Theorem \ref{th.KP.3} follows the next result.

\bco \label{cor.KP.4.2}
Under conditions of Theorem \ref{th.KP.3} we obtain
\beam \label{eq.KP.51}
\psi_{and}(x,\,y,\,n) \sim \sum_{k=1}^n \sum_{l=1}^n \PP\left[  X_k>x\,,\; Y_l>y \right]\,,
\eeam
as $x\wedge y \to \infty$.
\eco

\bre
From relation \eqref{eq.KP.51} and the definitions for  $T_{max}$ and $T_{and}$ we can easily observe that $\psi_{max}(x,\,y,\,n) \leq \psi_{and}(x,\,y,\,n)$, for any $x,\,y >0$ and any $n \in \bbn$. Thus, for  $\psi_{max}(x,\,y,\,n)$ we find the asymptotic upper bound
\beao
\psi_{max}(x,\,y,\,n) \lesssim  \sum_{k=1}^n \sum_{l=1}^n \PP\left[  X_k>x\,,\; Y_l>y \right]\,,
\eeao 
as $x\wedge y \to \infty$, for  any $n \in \bbn$.
\ere

\section{Scalar product} \label{sec.KP.5}

Now we examine the closure property of scalar product in $\mathcal{L}^{(2)}$, $\mathcal{D}^{(2)}$, and in their intersection. Later we check the same for random sums in two dimensions.

The scalar product has the following tail
\beam \label{eq.KP.5.2}
{\bf \bH} (x,\,y):=\PP[\Theta\,X > x,\,\Theta\,Y>y]\,.
\eeam
Here we set $\Theta$ to be a non-negative random variable with distribution $B$, such that $B(0-)=0$ and $B(0)<1$. We assume also that $\Theta$ is independent of $(X,\,Y)$.  These products in relation \eqref{eq.KP.5.2} have many applications in actuarial mathematics, in risk management, and stochastic fields. Next, we use an assumption from \cite{konstantinides:passalidis:2024}.

\begin{assumption} \label{ass.KP.A}
Let us suppose that it holds $\bB[c\,(x\wedge y)] = o\left(\PP\left[\Theta\,X>x\,,\;\Theta\,Y>y \right] \right)=:o\left[{\bf \bH} (x,\,y) \right]$, as $x\wedge y \to \infty$, for any $c>0$.
\end{assumption}

\bre \label{rem.KP.5.1}
From Assumption \ref{ass.KP.A} is  implied that
\beam \label{eq.KP.5.2b}
\bB(c\,x) = o\left(\PP\left[\Theta\,X>x\right] \right)\,,
\eeam
as $\xto$, for any $c>0$, and similarly 
\beam \label{eq.KP.5.3}
\bB(c\,y) = o\left(\PP\left[\Theta\,Y>y\right] \right)\,,
\eeam
as $\yto$ for any $c>0$. This condition is well-known, see for example \cite{tang:2006}. Further, we can see that Assumption  \ref{ass.KP.A} holds immediately when the distribution $B$ has support bounded from above.
\ere

The next lemma helps our argumentation and presents a multivariate extension of \cite[Lem. 3.2]{tang:2006}, providing the existence of an auxiliary function. Similar paper on auxiliary functions we refer to \cite{zhou:wang:wang:2012}.   

\ble \label{lem.KP.5.1}
For two distributions $B$ and ${\bf H} $, with $\bB(x)>0$, ${\bf \bH}(x,\,y)>0$ for any $x,\,y > 0$, Assumption  \ref{ass.KP.A} holds if and only if there exists a function $b\;:\;[0,\,\infty) \to (0,\,\infty)$, such that
\begin{enumerate}

\item
$b(x) \to \infty$, as $\xto$,

\item
$b(x)=o(x)$, as $\xto$,

\item
$\bB[b(x\wedge y)]=o[{\bf \bH}(x,\,y)]$, as $x\wedge y \to \infty$. 
\end{enumerate}
\ele

\pr~

($\Leftarrow$). The existence of such an auxiliary function easily implies Assumption \ref{ass.KP.A}, for example we consider the function $x\wedge y/n$.

($\Rightarrow$). Let suppose that Assumption \ref{ass.KP.A} is satisfied. Then we obtain
\beao
\lim_{x\wedge y \to \infty} \dfrac{\bB((x\wedge y) / n)}{{\bf \bH}(x,\,y)}=0\,.
\eeao
Let an increasing sequence of positive numbers $\{\lambda_n\,,\;n \in \bbn\}$ with $\lambda_{n+1}> (n+1)\,\lambda_n$, for any $n \in \bbn$, such that for any $x\wedge y \geq \lambda_n$ we have
\beao
\dfrac{\bB(x\wedge y/n)}{{\bf \bH}(x,\,y)} \leq \dfrac 1n
\,.
\eeao
Therefore, the points (1), (2) and (3) are satisfied with
\beao
b(x\wedge y):=\sup_{0\leq k \leq x\wedge y} z(k)\,, \qquad \qquad z(x\wedge y) = \sum_{n=1}^{\infty} \dfrac{x\wedge y}n\,{\bf 1}_{\{\lambda_n\leq x\wedge y \leq \lambda_{n+1} \}}\,,
\eeao
which completes the proof.
~\halmos

\bre
We have to mention that in case of distribution $B$ with support bounded from above, the existence of function $b(\cdot)$ follows immediately.
\ere

Now we study the closedness of class $\mathcal{D}^{(2)}$ under the scalar product.

\bth \label{th.KP.5.2}
Let $(X,\,Y)$ be random vector and  $\Theta$ be random variable, with tail distribution ${\bf \bF_1}(x\,y)=\PP[X > x,\,Y>y]$ and $B$ respectively, and assume $B(0-)=0$,$B(0)<1$. If  $\Theta$ and $(X,\,Y)$ are independent, Assumption  \ref{ass.KP.A} holds and $(F,G) \in \mathcal{D}^{(2)}$, then ${\bf H}(x,\,y) \in \mathcal{D}^{(2)}$.
\ethe

\pr~
Initially, we get from $F,\,G \in \mathcal{D}$ and by \cite[Th. 3.3 (i)]{cline:samorodnitsky:1994} or \cite[Prop. 5.4(i)]{leipus:siaulys:konstantinides:2023} that the products $\Theta\,X$ and $\Theta\,Y$ follow distributions from $\mathcal{D}$. From Assumption  \ref{ass.KP.A} we  obtain that for any ${\bf b} \in (0,\,1)^n$
\beam \label{eq.KP.5.9}
&&\limsup_{x\wedge y \to \infty}\dfrac {{\bf \overline{H}_b}(x,\,y)}{{\bf \overline{H}}(x,\,y)}=\limsup_{x\wedge y \to \infty}\dfrac{\PP[\Theta \,X>b_1\,x,\,\,\Theta \,Y>b_2\,y]}{\PP[\Theta \,X>x\,,\;\Theta \,Y>y]}=\limsup_{x\wedge y \to \infty}\\[4mm] \notag
&&\dfrac{\left(\int_0^{b(x\wedge y)}+\int_{b(x\wedge y)}^\infty \right)\PP\left[X>\dfrac{b_1\,x}{s}\,,\;Y>\dfrac{b_2\,y}s\right]\,B(ds)}{\PP[\Theta\, X>x\,,\;\Theta\, Y>y]}=:\limsup_{x\wedge y \to \infty}\dfrac{I_1+I_2}{\PP[\Theta\, X>x,\Theta\, Y>y]}\,.
\eeam
Further we calculate
\beao
I_2 \leq \int_{b(x\wedge y)}^\infty \,B(ds)=\bB[b(x\wedge y)] = o\left[{\bf \overline{H}}(x,\,y)\right]\,,
\eeao
as $x\wedge y \to \infty$, due to Assumption  \ref{ass.KP.A}. Hence, taking into account also relation \eqref{eq.KP.5.9} we find
\beao
\limsup_{x\wedge y \to \infty}\dfrac {{\bf \overline{H}_b}(x,\,y)}{{\bf \overline{H}}(x,\,y)}&\leq&\limsup_{x\wedge y \to \infty}\dfrac{\int_0^{b(x\wedge y)}\PP\left[X>\dfrac{b_1\,x}{s}\,,\;Y>\dfrac{b_2\,y}s\right]\,B(ds)}{\int_0^{b(x\wedge y)}\PP\left[X>\dfrac{x}{s}\,,\;Y>\dfrac{y}s\right]\,B(ds)}\\[2mm] 
&\leq&\limsup_{x\wedge y \to \infty} \sup_{0< s \leq b(x\wedge y)}\dfrac{\PP\left[X>b_1\,x/s,\;Y>b_2\,y/s\right]}{\PP\left[X_1>x/s,\;Y>y/s\right]}\\[2mm] 
&\leq&\limsup_{x\wedge y \to \infty}\dfrac{\PP\left[X>b_1\,x,\;Y>b_2\,y\right]}{\PP\left[X> x,\;Y> y\right]}<\infty\,,
\eeao
where in the last step we used the condition $(F,G) \in \mathcal{D}^{(2)}$. So we get ${\bf H}(x,\,y) \in \mathcal{D}^{(2)}$.~\halmos

Let us observe, that if $\Theta$ has upper bounded support, the proof of Theorem \ref{th.KP.5.2} (as also of  Theorem \ref{th.KP.5.1}) is implied by similar manipulations, replacing $b\,(x\wedge y)$ by the right endpoint of the distribution $B$. Next, we provide an analogue for class $\mathcal{L}^{(2)}$.

\bth \label{th.KP.5.1}
Let $(X,\,Y)$ be a random vector and $\theta$ be a non-negative random variable, with distributions ${\bf F}$, $B$ respectively, under condition $B(0) <1$. If $(X,\,Y)$, $\Theta$ are independent, Assumption \ref{ass.KP.A} holds, and $(F,\,G) \in \mathcal{L}^{(2)}$, then ${\bf H}(x,\,y) \in \mathcal{L}^{(2)}$.
\ethe

\pr~
From the fact that $(X,\,Y)$ is independent of $\Theta$, $F,\,G \in \mathcal{L}$ and relations \eqref{eq.KP.5.2b} and \eqref{eq.KP.5.3}, using \cite[Th 2.2 (iii)]{cline:samorodnitsky:1994}, we find that distributions of $\Theta\,X$ and $\Theta\,Y$ belong to $\mathcal{L}$. Let ${\bf a}=(a_1,\,a_2) >(0,\,0)$. Then we easily obtain
\beam \label{eq.KP.5.4}
\liminf_{x\wedge y \to \infty}\dfrac{{\bf \bH} (x-a_1,\,y-a_2)}{{\bf \bH} (x,\,y) }=\lim_{x\wedge y \to \infty}\dfrac {\PP[\Theta\,X > x-a_1,\,\Theta\,Y>y-a_2]}{\PP[\Theta\,X >x,\;\Theta\,Y>y]} \geq 1\,.
\eeam
Next, we show the opposite asymptotic inequality. Using Assumption \ref{ass.KP.A} we obtain
\beam \label{eq.KP.5.5}
&&\limsup_{x\wedge y \to \infty}\dfrac{{\bf \bH} (x-a_1,\,y-a_2)}{{\bf \bH} (x,\,y) }\\[2mm] \notag
&&=\lim_{x\wedge y \to \infty}\dfrac {1}{{\bf \bH}(x,\,y)} \,\left(\int_0^{b(x\wedge y)} + \int_{b(x\wedge y)}^\infty \right)\PP\left[ X>\dfrac{x -a_1}{s},\;Y>\dfrac{y-a_2}{s} \right]B(ds)\\[2mm] \notag
&&=:\lim_{x\wedge y \to \infty}\dfrac {I_1(x,\,y)  +I_2(x,\,y) }{{\bf \bH}(x,\,y) } \,.
\eeam
Thus, by Assumption \ref{ass.KP.A} we find
\beao
I_2(x,\,y)  =  \int_{b(x\wedge y)}^\infty \,\PP\left[ X>\dfrac{x -a_1}{s}\,,\;Y>\dfrac{y-a_2}{s} \right]\,B(ds) \leq \bB[b(x\wedge y)] = o\left[ {\bf \bH} (x,\,y) \right]\,,
\eeao
hence, 
\beao
\dfrac{I_2(x,\,y)}{{\bf \bH} (x,\,y)} = o(1)\,,
\eeao 
as $x\wedge y \to \infty$. As a consequence, taking into account also \eqref{eq.KP.5.5} we get
\beao
\limsup_{x\wedge y \to \infty}\dfrac{{\bf \bH_1} (x-a_1,\,y-a_2)}{{\bf \bH} (x,\,y) }&=&\limsup_{x\wedge y \to \infty}\int_0^{b(x\wedge y)} \PP\left[ X>\dfrac{x -a_1}{s},\;Y>\dfrac{y -a_2}{s} \right]\dfrac {B(ds)}{{\bf \bH} (x,\,y)} \\[2mm]
&\leq& \limsup_{x\wedge y \to \infty}\dfrac {\int_0^{b(x\wedge y)} \PP\left[ X>\dfrac{x -a_1}{s},\;Y>\dfrac{y -a_2}{s} \right]B(ds)}{\int_0^{b(x\wedge y)} \PP\left[ X>\dfrac{ x }{s},\;Y>\dfrac{y }{s} \right]B(ds)}\\[2mm]
&\leq& \limsup_{x\wedge y \to \infty} \sup_{0<s \leq b(x\wedge y)}\dfrac {\PP\left[ X>\dfrac{x -a_1}{s},\;Y>\dfrac{y-a_2}{s} \right]}{\PP\left[ X>\dfrac{ x }{s},\;Y>\dfrac{y }{s} \right]} \\[2mm]
&=& \limsup_{x\wedge y \to \infty} \dfrac {\PP\left[ X>x -a_1,\;Y>y -a_2 \right]}{\PP\left[ X>x\,,\;Y>y \right]}=1 \,.
\eeao
where in the last step we consider the fact that  $(F,G) \in \mathcal{L}^{(2)}$. So we have 
\beam \label{eq.KP.5.6}
\limsup_{x\wedge y \to \infty}\dfrac{{\bf \bH} (x-a_1,\,y-a_2)}{{\bf \bH} (x,\,y) }\leq 1\,.
\eeam
From relations \eqref{eq.KP.5.4} and \eqref{eq.KP.5.6} we conclude  ${\bf H} (x,\,y) \in \mathcal{L}^{(2)}$.
~\halmos

The next statement stems from a combination of previous results. 

\bco \label{cor.KP.5.1}
Let  $(X,\,Y)$ be a random vector and $\Theta$ be a non-negative random variable with distributions $(F,G)$, $B$ respectively, under condition $B(0) <1$. If $(X,\,Y)$ and $\Theta$ are independent, with 
 $(X,\,Y)\in(\mathcal{D}\cap\mathcal{L})^{(2)}$ and satisfy the Assumption  \ref{ass.KP.A}, then ${\bf H} (x,\,y) \in (\mathcal{D}\cap\mathcal{L})^{(2)}$.
\eco

\pr~
This follows directly from Theorem \ref{th.KP.5.2} and Theorem \ref{th.KP.5.1}. 
~\halmos

\section{Randomly weighted sums} \label{sec.KP.6}

Finally, we extend Theorem \ref{th.KP.3} into weighted sums. The first kind of weighted sums takes the form 
\beao
S_n(\Theta)=\sum_{k=1}^n \Theta\,X_k\,, \qquad T_n(\Theta)=\sum_{l=1}^n \Theta\,Y_l\,.
\eeao 
These quantities have the same discount factor $\Theta$, hence the $(X_k,\,Y_l)$, for $k,\,l=1,\,\ldots,\,n$, are the losses or gains of the two lines of business during the $k$-th period. If the $(x,\,y)$ represents the two initial capitals respectively, then  the ruin probability in this model comes in the form
\beam \label{eq.KP.6.2}
\psi_{and}(x,\,y,\,n):=\PP\left[\bigvee_{i=1}^n S_i(\Theta)>x\,,\;\bigvee_{j=1}^n T_j(\Theta)>y \right]\,.
\eeam
The ruin probability plays a significant role in risk theory. For example we refer to, \cite{li:tang:2015}, \cite{yang:konstantinides:2015} and \cite{cheng:2021}, \cite{ji:wang:yan:cheng:2023} for discrete-time or continuous-time models respectively.

The next result is based  on Theorem \ref{th.KP.3} and Corollary \ref{cor.KP.5.1}. We have to notice that there exists the asymptotic behavior of the ruin probability in \eqref{eq.KP.6.2} as well.

\bco \label{cor.KP.6.1}
Let $X_1,\,\ldots,\,X_n,\,Y_1,\,\ldots,\,Y_n$ be random variables with the following distributions $F_1,\,\ldots,\,F_n,\,G_1,\,\ldots,\,G_n$ respectively form class $\mathcal{D}\cap\mathcal{L}$ and they satisfy the GTAI dependence structure. We assume that $\Theta$ represents a non-negative random variable with upper bounded support, it is independent of $X_1,\,\ldots,\,X_n,\,Y_1,\,\ldots,\,Y_n$ and $(X_k,\,Y_l) \in (\mathcal{D}\cap\mathcal{L})^{(2)}$, for $k,\,l=1,\,\ldots,\,n$. If further $X_1,\,\ldots,\,X_n$ are TAI and $Y_1,\,\ldots,\,Y_n$ are TAI then the following asymptotic relation is true
\beam \label{eq.KP.6.4}
&&\PP\left[ S_n(\Theta)>x\,,\; T_n(\Theta)>y \right] \sim \PP\left[\bigvee_{i=1}^n S_i(\Theta)>x\,,\;\bigvee_{j=1}^n T_j(\Theta)>y \right] \\[2mm]  \notag
&&\sim\PP\left[\bigvee_{k=1}^n \Theta\,X_k>x\,,\;\bigvee_{l=1}^n \Theta\,Y_l>y \right]  \sim  \sum_{k=1}^n  \sum_{l=1}^n \PP\left[ \Theta\,X_k>x\,,\;\Theta\,Y_l>y \right] \,,
\eeam
as $x\wedge y \to \infty$.
\eco

\pr~
We start from \cite[Lem. 3.1]{konstantinides:passalidis:2023} and because of the upper-bound of $\Theta$, we obtain that the products $\Theta\,X_1,\,\ldots,\,\Theta\,X_n,\,\Theta\,Y_1,\,\ldots,\,\Theta\,Y_n$ are GTAI. Now we can apply \cite[Th. 2.2(iii), Th.3.3(ii)]{cline:samorodnitsky:1994}, to find $\Theta\,X_k \in \mathcal{D}\cap\mathcal{L}$, and  $\Theta\,Y_l \in \mathcal{D}\cap\mathcal{L}$ for any  $k=1,\,\ldots,\,n$ and  $l=1,\,\ldots,\,n$. Because of the closedness of class $\mathcal{D}$ we using Theorem 2.2 of \cite{li:2013} $\Theta_1\,X_1,\,\ldots,\,\Theta_n\,X_n$ are TAI and $\Delta_1\,Y_1,\,\ldots,\,\Delta_n\,Y_n$ are TAI.

Next, since $\Theta$ is bounded from above, so Assumption \ref{ass.KP.A} is fulfilled for $\Theta$, $(X_k,\,Y_l)$, for any $k,\,l=1,\,\ldots,\,n$, in order to obtain $(\Theta\,X_k,\,\Theta\,Y_l) \in \mathcal{D}^{(2)}$ for any  $k=1,\,\ldots,\,n$ and  $l=1,\,\ldots,\,n$, it is enough to apply Theorem \ref{th.KP.5.2}, and similarly, by Theorem \ref{th.KP.5.1}, we find $(\Theta\,X_k,\,\Theta\,Y_l) \in \mathcal{L}^{(2)}$ for any  $k=1,\,\ldots,\,n$ and  $l=1,\,\ldots,\,n$. Therefore, the $(\Theta\,X_k,\,\Theta\,Y_l) \in (\mathcal{D}\cap\mathcal{L})^{(2)}$  and the $\Theta\,X_1,\,\ldots,\,\Theta\,X_n$, $\Theta\,Y_1,\,\ldots,\,\Theta\,Y_n$ are GTAI. Finally, applying Theorem \ref{th.KP.3}, we conclude \eqref{eq.KP.6.4}.
~\halmos

Now we need some preliminary results. Several times before proving that the convolution product satisfies $H \in \mathcal{B}$, with $\mathcal{B}$ some distribution class, we need to prove that $H_{\vep}(x):=\PP[(\Theta \vee \vep)\,X \leq x]$ belongs to this class $\mathcal{B}$ for any $\vep>0$. Following the approach in \cite{cline:samorodnitsky:1994} we show that for some constant $\delta>0$, if $H_{\vep} \in \mathcal{L}^{(2)}$, for any $\vep \in (0,\,\delta)$, then $H \in \mathcal{L}^{(2)}$. However the next results, deserves theoretical attention by its own merit. 

From here until to the end of paper we assume that $X,Y$ are non-negative random variables.

\ble \label{lem.KP.6.1}
If for some constant vector ${\bf \delta}=(\delta_1,\,\delta_2) > (0,\,0)$, for any $\vep_1 \in (0,\,\delta_1)$  and for any $\vep_2 \in (0,\,\delta_2)$, holds $((\Theta\vee \vep_1)\,X,\,(\Delta\vee \vep_2)\,Y) \in\mathcal{L}^{(2)}$, with $X,Y,\Theta,\Delta$ non-negative random variables, then we conclude that $(\Theta\,X,\,\Delta\,Y) \in\mathcal{L}^{(2)}$. 
\ele

\pr~
Keeping in mind that $((\Theta\vee \vep_1)\,X,\,(\Delta\vee \vep_2)\,Y) \in\mathcal{L}^{(2)}$, we start by \cite[th. 2.2(i)]{cline:samorodnitsky:1994} to establish that due to $(\Theta\vee \vep_1)\,X \in\mathcal{L}\,,\; (\Delta\vee \vep_2)\,Y \in\mathcal{L}$, we get  $\Theta\,X \in\mathcal{L}$, and $\Delta\,Y\in\mathcal{L}$. Next we check the second property of class $\mathcal{L}^{(2)}$. Let $(a_1,\,a_2) > (0,\,0)$, then
\beam \label{eq.KP.6.5} 
\liminf_{x\wedge y \to \infty}\dfrac{\PP\left[\Theta\,X >x-a_1\,,\;\Delta\,Y>y -a_2 \right]}{\PP\left[\Theta\,X >x\,,\;\Delta\,Y>y \right]} \geq 1 \,,
\eeam
Next, for any $(\vep_1,\,\vep_2) > (0,\,0)$, we find
\beao
&&\PP[(\Theta\vee \vep_1)\,X>x\,,\;(\Delta\vee \vep_2)\,Y>y] \geq \PP[\Theta\,X>x\,,\;\Delta\,Y>y]\geq \PP[\Theta\,X>x\,,\;\Theta>\vep_1\,,\\[2mm]
&&\;\Delta\,Y>y]=\PP[(\Theta\vee \vep_1)\,X>x\,,\;\Delta\,Y>y]-\PP[\Theta \leq \vep_1]\,\PP[X\,\vep_1>x,\,\Delta\,Y>y]\\[2mm]
&&\geq \PP[\Theta > \vep_1]\,\PP[(\Theta\vee \vep_1)\,X>x,\,\Delta\,Y>y,\,\Delta>\vep_2]= \PP[\Theta > \vep_1]\,\\[2mm]
&&\times \left(\PP[(\Theta\vee \vep_1)\,X>x,\,(\Delta\vee \vep_2)\,Y>y]-\PP[(\Theta\vee \vep_1)\,X>x,\, \vep_2\,Y>y] \,\PP[\Delta \leq \vep_2]\right)\\[2mm]
&&\geq \PP[\Theta > \vep_1]\,\PP[\Delta > \vep_2]\,\PP[(\Theta\vee \vep_1)\,X>x,\,(\Delta\vee \vep_2)\,Y>y]\,,
\eeao
hence we conclude
\beam \label{eq.KP.6.6}
&&\PP[(\Theta\vee \vep_1)\,X>x\,,\;(\Delta\vee \vep_2)\,Y>y] \geq \PP[\Theta\,X>x\,,\;\Delta\,Y>y]\\[2mm] \notag
&&\geq \PP[\Theta > \vep_1]\,\PP[\Delta > \vep_2]\,\PP[(\Theta\vee \vep_1)\,X>x,\,(\Delta\vee \vep_2)\,Y>y]\,.
\eeam
Therefore, using \eqref{eq.KP.6.6} and due to properties of $\mathcal{L}^{(2)}$, for $((\Theta\vee \vep_1)\,X,\,(\Delta\vee \vep_2)\,Y)$ we obtain
\beao
&&\limsup_{x\wedge y \to \infty}\dfrac{\PP\left[\Theta\,X >x-a_1\,,\;\Delta\,Y>y -a_2 \right]}{\PP\left[\Theta\,X >x\,,\;\Delta\,Y>y \right]} \leq \limsup_{x\wedge y \to \infty}\\[2mm]
&&\dfrac{\PP\left[(\Theta\vee \vep_1)\,X >x-a_1\,,\;(\Delta\vee \vep_2)\,Y>y -a_2 \right]}{\PP[\Theta > \vep_1] \PP[\Delta > \vep_2]\,\PP\left[(\Theta\vee \vep_1) X >x ,\,(\Delta\vee \vep_2) Y>y\right]}=\dfrac{1}{\PP[\Theta > \vep_1] \PP[\Delta > \vep_2]}\,,
\eeao
and leaving $\vep_1$ and $\vep_2$ to tend to zero, we get
\beao
\limsup_{x\wedge y \to \infty}\dfrac{\PP\left[\Theta\,X >x-a_1\,,\;\Delta\,Y>y -a_2 \right]}{\PP\left[\Theta\,X >x\,,\;\Delta\,Y>y \right]} \leq 1\,,
\eeao
hence, from \eqref{eq.KP.6.5} and from last inequality we reach to $(\Theta\,X,\,\Delta\,Y) \in\mathcal{L}^{(2)}$.
~\halmos

\ble \label{lem.KP.6.2}
Let $X$ and $Y$ be non-negative random variables, with $(X,\,Y) \in \mathcal{L}^{(2)}$ and $\Theta$ and $\Delta$ be non-negative, non-degenerated to zero random variables, independent of $(X,\,Y)$. We assume that 
\beam \label{eq.KP.6.8} 
&&\PP\left[\Theta >x\right]=o(\PP\left[\Theta\,X >c_1\,x\,,\;\Delta\,Y>c_2\,y \right])= \PP\left[\Delta >y\right]\,,
\eeam
as $x \wedge y \to \infty$, for any $c_1,\,c_2 >0$. Then $(\Theta\,X,\,\Delta\,Y) \in\mathcal{L}^{(2)}$.
\ele

\pr~
From \eqref{eq.KP.6.8} we obtain
\beam \label{eq.KP.6.8*}
\dfrac{\PP\left[\Theta >x \right]}{\PP\left[\Theta \,X >c_1\,x\right]} \leq \dfrac{\PP\left[\Theta >x \right]}{\PP\left[\Theta \,X >c_1\,x\,,\;\Delta\,Y>c_2\,y\right]} \longrightarrow 0\,,
\eeam
as $x\wedge y \to \infty$, and similarly we find $\PP\left[\Delta >y\right]=o(\PP\left[\Delta\,Y>c_2\,y \right])$, as $x\wedge y \to \infty$, for any $c_1,\,c_2 >0$. Hence, by \cite[Th. 2.2]{cline:samorodnitsky:1994} we find $\Theta \,X  \in\mathcal{L}$ and $\Delta\,Y \in\mathcal{L}$. Next, we show the second property of $(\Theta\,X,\,\Delta\,Y) \in\mathcal{L}^{(2)}$. Indeed, from Lemma \ref{lem.KP.6.1} we see that it is enough to show this for any $\Theta \geq \vep_1$ and $\Delta \geq \vep_2$ almost surely for any $\vep_1,\,\vep_2 >0$. Let consider some $a_1,\,a_2 >0$ and some $k_1,\,k_2,\,k >0$, such that for a large enough $x_0 \geq 0$, it holds
\beam \label{eq.KP.6.9} 
\PP\left[X >x-\dfrac{a_1}{\vep_1}\,,\;Y>y -\dfrac{a_2}{\vep_2} \right]\leq (1+k)\,\PP\left[X >x\,,\;Y> y \right]\,,
\eeam
for any $x \wedge y \geq x_0$ and 
\beam \label{eq.KP.6.10} 
\PP\left[X >x-\dfrac{a_1}{\vep_1}\right]\leq (1+k_1)\,\PP\left[X >x \right]\,,\quad \PP\left[Y>y -\dfrac{a_2}{\vep_2} \right]\leq (1+k_2)\,\PP\left[Y> y \right]\,,
\eeam
for any $x \geq x_0$ and $y \geq x_0$, respectively. Then, for all $x \wedge y \geq x_0$, we obtain
\beam \label{eq.KP.6.11} 
&&\PP\left[\Theta\,X >x-a_1\,,\;\Delta\,Y>y -a_2 \right]=\left(\int_{\vep_1}^{x/x_0}+\int_{x/x_0}^\infty \right)\,\left(\int_{\vep_2}^{y/x_0}+\int_{y/x_0}^\infty \right)\,\\[2mm] \notag
&&\qquad \qquad \PP\left[X >\dfrac{x-a_1}{s}\,,\;Y>\dfrac{y-a_2}{t} \right] \,\PP[\Theta \in ds,\;\Delta \in dt]=:\sum_{m=1}^4 I_m(x,\,y)\,,
\eeam
where we find
\beao
I_4(x,\,y) = \int_{x/x_0}^\infty \,\int_{y/x_0}^\infty\PP\left[X >\dfrac{x-a_1}{s}\,,\;Y>\dfrac{y-a_2}{t} \right] \,\PP[\Theta \in ds,\;\Delta \in dt]\,,
\eeao
that gives 
\beam \label{eq.KP.6.12} 
I_4(x,\,y) \leq \PP\left[\Theta \geq \dfrac x{x_0}\,,\;\Delta  \geq \dfrac y{x_0}\right]\,.
\eeam 

Now we estimate $I_1(x,\,y)$
\beam  \label{eq.KP.6.13} \notag
I_1(x,\,y) &=& \int_{\vep_1}^{x/x_0}\,\int_{\vep_2}^{y/x_0} \PP\left[X >\dfrac{x-a_1}{s}\,,\;Y>\dfrac{y-a_2}{t} \right] \,\PP[\Theta \in ds,\;\Delta \in dt]\\[2mm]  \notag
&\leq& \int_{\vep_1}^{x/x_0}\,\int_{\vep_2}^{y/x_0}\,\PP\left[X >\dfrac{x}{s}-\dfrac{a_1}{\vep_1}\,,\;Y>\dfrac{y}{t}-\dfrac{a_2}{\vep_2} \right] \,\PP[\Theta \in ds,\;\Delta \in dt]\\[2mm] 
&\leq& (1+k)\,\int_{\vep_1}^{x/x_0}\,\int_{\vep_2}^{y/x_0}\,\PP\left[X >\dfrac{x}{s}\,,\;Y>\dfrac{y}{t} \right] \,\PP[\Theta \in ds,\;\Delta \in dt]\\[2mm]  \notag
&\leq& (1+k)\,\PP\left[\Theta \,X >x\,,\;\Delta\,Y>y \right] \,,
\eeam
thus we get $I_1(x,\,y) \leq  (1+k)\,\PP\left[\Theta \,X> x\,,\;\Delta\,Y >y\right]$, which follows from \eqref{eq.KP.6.9}.

Next we consider $I_2(x,\,y)$
\beam  \label{eq.KP.6.14} \notag
I_2(x,\,y) &=& \int_{\vep_1}^{x/x_0}\,\int_{y/x_0}^\infty \PP\left[X >\dfrac{x-a_1}{s}\,,\;Y>\dfrac{y-a_2}{t} \right] \,\PP[\Theta \in ds,\;\Delta \in dt]\\[2mm] \notag 
&\leq& \int_{\vep_1}^{x/x_0}\,\PP\left[X >\dfrac{x}{s}-\dfrac{a_1}{\vep_1} \right] \,\PP\left[\Theta \in ds,\;\Delta > \dfrac y{x_0} \right]\\[2mm] 
&\leq& (1+k_1)\, \int_{\vep_1}^{x/x_0}\,\PP\left[X >\dfrac{x}{s} \right] \,\PP\left[\Theta \in ds,\;\Delta > \dfrac y{x_0} \right]\\[2mm]  \notag
&\leq& (1+k_1)\,\PP\left[\Theta \,X >x\,,\;\Delta > \dfrac y{x_0} \right] \leq (1+k_1)\,\PP\left[\Delta > \dfrac y{x_0} \right]\,,
\eeam
that means $I_2(x,\,y) \leq  (1+k_1)\,\PP\left[\Delta>y/{x_0} \right]$, where in the pre-last step we use the first relation in  \eqref{eq.KP.6.10}. 

For $I_3(x,\,y)$, we use the second relation in  \eqref{eq.KP.6.10} and due to symmetry with respect to \eqref{eq.KP.6.14}, we find
\beam \label{eq.KP.6.15} \notag
I_3(x,\,y) &=&\int_{x/x_0}^{\infty} \int_{\vep_2}^{\infty}\, \PP\left[X >\dfrac{x-a_1}{s}\,,\;Y>\dfrac{y-a_2}{t} \right] \,\PP[\Theta \in ds,\;\Delta \in dt]\\[2mm] 
&\leq&  (1+k_2)\,\PP\left[\Theta> \dfrac x{x_0} \right]\,.
\eeam

Therefore putting the estimation from \eqref{eq.KP.6.12} - \eqref{eq.KP.6.15} into \eqref{eq.KP.6.11} we conclude
\beao
&&\PP\left[\Theta\,X >x-a_1\,,\;\Delta\,Y>y -a_2 \right] \leq\PP\left[\Theta> \dfrac x{x_0}\,,\;\Delta> \dfrac y{x_0}  \right] \\[2mm] \notag
&&\qquad  + (1+k)\,\PP\left[\Theta\,X>x\,,\;\Delta\,Y>y \right]+ (1+k_2)\,\PP\left[\Theta> \dfrac x{x_0} \right]+  (1+k_1)\,\PP\left[\Delta> \dfrac y{x_0} \right]\,,
\eeao
Now, because of \eqref{eq.KP.6.8} and the relation 
\beao
\dfrac{\PP\left[\Theta >x\,,\; \Delta> y\right]}{\PP\left[\Theta \,X >x\,,\;\Delta\,Y>y\right]} \leq \dfrac{\PP\left[\Theta >x \right]}{\PP\left[\Theta \,X >x\,,\;\Delta\,Y>y\right]} \longrightarrow 0\,,
\eeao
as $x \wedge y \to \infty$, we find
\beao 
\lim_{x\wedge y \to \infty}\dfrac{\PP\left[\Theta\,X >x-a_1\,,\;\Delta\,Y>y -a_2 \right]}{\PP\left[\Theta\,X >x\,,\;\Delta\,Y>y \right]} \leq 1+k\,.
\eeao
By these inequalities and relation \eqref{eq.KP.6.8*}, in combination of the arbitrary choice of $k$ and relation \eqref{eq.KP.6.5} we have $(\Theta\,X,\,\Delta\,Y) \in\mathcal{L}^{(2)}$.  
~\halmos

Next, we consider the asymptotic joint-tail behavior of discounted aggregate claims in a two-dimensional risk model on discrete-time, where the vector $(X_k,\,Y_k)$ represents losses in two lines of business at the $k$-th period, while the $(\Theta_k,\,\Delta_k)$ represent the discount factors of these two lines of business respectively. In this risk model we study only the aggregate claims, and we accept that the $\Theta_1,\,\ldots,\,\Theta_n,\,\Delta_1,\,\ldots,\,\Delta_m$ are independent of claims $X_1,\,\ldots,\,X_n,\,Y_1,\,\ldots,\,Y_m$. For further reading on risk models with dependence among the discount factors and main claims see in \cite{chen:2011}, \cite{chen:2017}, \cite{yang:gao:li:2016}, but only in one dimension. Namely we have the sums:
\beao
S_{n}^{\Theta}:=\sum_{k=1}^{n}\Theta_k X_k\,,\qquad T_{n}^{\Delta}:=\sum_{l=1}^{n}\Delta_l Y_l\,.
\eeao

\begin{assumption} \label{ass.KP.C}
There exist constants  $0<\xi_k \leq \delta_k$ such that hold $\xi_k\leq \Theta_k \leq \delta_k$ almost surely, for any $k=1,\,\ldots,\,n$ and there exist constants $0<\gamma_l \leq \zeta_l$ such that hold $\gamma_l\leq \Delta_l \leq \zeta_l$ almost surely, for any $l=1,\,\ldots,\,n$.
\end{assumption}

\bth \label{th.KP.6.1}
Let $X_1,\,\ldots,\,X_n,\,Y_1,\,\ldots,\,Y_n$ be random variables with the following distributions $F_1,\,\ldots,\,F_n,\,G_1,\,\ldots,\,G_n$ respectively form class $\mathcal{D}\cap\mathcal{L}$ and they satisfy the GTAI dependence structure, with $(X_k,\,Y_l) \in (\mathcal{D}\cap\mathcal{L})^{(2)}$ for any $k=1,\,\ldots,\,n$ and $l=1,\,\ldots,\,n$. We suppose that the random discount factors  $\Theta_1,\,\ldots,\,\Theta_n,\,\Delta_1,\,\ldots,\,\Delta_n$ satisfy Assumption \ref{ass.KP.C} and are independent of $X_1,\,\ldots,\,X_n,\,Y_1,\,\ldots,\,Y_n$. Then the  products $\Theta_1\,X_1,\,\ldots,\,\Theta_n\,X_n$, $\Delta_1\,Y_1,\,\ldots,\,\Delta_n\,Y_n$, are GTAI with $(\Theta_k\,X_k,\,\Delta_l\,Y_l) \in (\mathcal{D}\cap\mathcal{L})^{(2)}$. If further $X_1,\,\ldots,\,X_n$ are TAI and $Y_1,\,\ldots,\,Y_n$ are TAI then holds the asymptotic relations
\beam \label{eq.KP.6.20} 
&&\PP\left[S_n^{\Theta}>x\,,\;T_n^{\Delta}>y\right] \sim\PP\left[\bigvee_{i=1}^n S_i^{\Theta}>x\,,\;\bigvee_{j=1}^n T_j^{\Delta}>y\right] \\[2mm] \notag
&&\qquad \qquad \sim \PP\left[\bigvee_{k=1}^n \Theta_k\,X_k>x\,,\;\bigvee_{l=1}^n \Delta_l\,Y_l>y\right] \sim \sum_{k=1}^n  \sum_{l=1}^n \PP\left[ \Theta\,X_k>x\,,\;\Delta\,Y_l>y \right] 
\,,
\eeam
as $x\wedge y \to \infty$.
\ethe

\pr~
Taking into account the upper bound for discount factors $\Theta_1,\,\ldots,\,\Theta_n,\,\Delta_1,\,\ldots,\,\Delta_n$  and their independence from  $X_1,\,\ldots,\,X_n,\,Y_1,\,\ldots,\,Y_n$, we apply \cite[Lem. 1]{konstantinides:passalidis:2023} to find that the products $\Theta_1\,X_1,\,\ldots,\,\Theta_n\,X_n,\,\Delta_1\,Y_1,\,\ldots,\,\Delta_n\,Y_n$ are GTAI. Now by \cite[Th. 3.3(i)]{cline:samorodnitsky:1994} we get $\Theta_k\,X_k \in \mathcal{D}\cap\mathcal{L}$ and $\Delta_l\,Y_l \in \mathcal{D}\cap\mathcal{L}$, for any $k=1,\,\ldots,\,n$ and for any $l=1,\,\ldots,\,n$.  As a result by class $\mathcal{D}$ we using Theorem 2.2 of \cite{li:2013} $\Theta_1\,X_1,\,\ldots,\,\Theta_n\,X_n$ are TAI and $\Delta_1\,Y_1,\,\ldots,\,\Delta_n\,Y_n$ are TAI.

Next, we check if $(\Theta_k\,X_k,\,\Delta_l\,Y_l) \in (\mathcal{D}\cap\mathcal{L})^{(2)}$ for any $k=1,\,\ldots,\,n$ and $l=1,\,\ldots,\,n$. Let ${\bf b}=(b_1,\,b_2) \in (0,\,1)^2$, then
\beao
\limsup_{x\wedge y \to \infty}\dfrac{\PP\left[\Theta_k\,X_k >b_1\,x,\;\Delta_l\,Y_l>b_2\,y  \right]}{\PP\left[\Theta_k\,X_k >x,\,\Delta_l\,Y_l>y \right]}\leq\limsup_{x\wedge y \to \infty}\dfrac{\PP\left[X_k >b_1\dfrac x{\delta_k},\,Y_l>b_2\dfrac y{\zeta_l}  \right]}{\PP\left[X_k >\dfrac x{\xi_k}\,,\;Y_l>\dfrac y{\gamma_l} \right]} <\infty,
\eeao
which follows from the inequalities 
\beao
\dfrac {b_1}{\delta_k} < \dfrac 1{\xi_k}\,, \qquad \dfrac {b_2}{\zeta_l} < \dfrac 1{\gamma_l}\,,
\eeao 
and the membership $(X_k,\,Y_l) \in (\mathcal{D}\cap\mathcal{L})^{(2)}$ for any $k=1,\,\ldots,\,n$ and $l=1,\,\ldots,\,n$. Next, we find the relation  $(\Theta_k\,X_k,\,\Delta_l\,Y_l) \in\mathcal{D}^{(2)}$ for any $k=1,\,\ldots,\,n$ and $l=1,\,\ldots,\,n$.

Now, noticing that relations \eqref{eq.KP.6.8} are satisfied because of Assumption  \ref{ass.KP.C}, we obtain directly from Lemma \ref{lem.KP.6.2} the inclusion $(\Theta_k\,X_k,\,\Delta_l\,Y_l) \in \mathcal{L}^{(2)}$ for any $k=1,\,\ldots,\,n$ and $l=1,\,\ldots,\,n$. Hence $(\Theta_k\,X_k,\,\Delta_l\,Y_l) \in (\mathcal{D}\cap\mathcal{L})^{(2)}$ for any $k=1,\,\ldots,\,n$ and $l=1,\,\ldots,\,n$ and by application of Theorem \ref{th.KP.3} for the products we conclude relation \eqref{eq.KP.6.20}.
~\halmos

\noindent \textbf{Acknowledgments.} 
We feel the pleasant duty to express deep gratitude to anonymous referees, who gave concise advices, that improved significantly the paper.

{\bf Competing Interest.} None.

{\bf Data availability.} There are no data used in this research.

{\bf Funding.} None.

\end{document}